\documentclass[12pt,a4paper]{article}
\usepackage{amssymb}

\newcommand{\plun}{\Delta_{p_1}}
\newcommand{\pldeux}{\Delta_{p_2}}
\newcommand{\Om}{\overline{\Omega}}
\newcommand{\R}{\mathbb{R}}
\newcommand{\io}{\int_{\Omega}}

\renewcommand{\epsilon}{\varepsilon}
\newcommand{\pl}{\Delta_{p}}
\newcommand{\N}{\mathbb{N}}
\newcommand{\espaces}{\R^+\times C^1(\Om)\times C^1(\Om)}
\renewcommand{\phi}{\varphi}

\newcommand{\cqfd}{\mbox{}\nolinebreak\hfill$\blacksquare$\medbreak\par}

\newtheorem{thm}{Theorem}[section]
\newtheorem{cor}[thm]{Corollary}
\newtheorem{lemma}[thm]{Lemma}

\newtheorem{remark}{Remark}[section]
\newtheorem{definition}{Definition}[section]
\newcommand{\proof}{\textbf{Proof: }}
\begin{document}
\title{
Symmetry and monotonicity results for  positive\\ solutions of p-Laplace
systems}
\author{C\'eline Azizieh}\date{}
\maketitle
\abstract{In this paper, we  extend to a system of
the type:
$$
\left\{\begin{array}{llll}
-\plun u= f(v)&\textrm{in }\Omega,& u>0\quad\textrm{in
}\Omega,&u=0\quad\textrm{on }\partial\Omega,\vspace*{0.2cm}\\
-\pldeux v= g(u)&\textrm{in }\Omega,& v>0\quad\textrm{in
}\Omega,&v=0\quad\textrm{on
}\partial\Omega,
\end{array}\right.
$$
where $\Omega\subset\R^N$ is
bounded, the monotonicity and symmetry results of
Damascelli and Pacella obtained in \cite{Da} in the case of a scalar
p-Laplace equation
with $1<p<2$. For this purpose, we use
the moving hyperplanes method and we suppose that
$f,g:\R\to\R^+$ are  increasing on $\R^+$
and
locally Lipschitz continuous
on $\R$ and $p_1,p_2\in(1,2)$ or $p_1\in(1,\infty),
 p_2=2$.}
 \mbox{}\\

\section{Introduction and statement of the main results}
Let $\Omega\subset\R^N$ be a bounded domain with $C^1$ boundary and
let
$f,g:\R\to\R^+$ be increasing on $\R^+$,
locally Lipschitz continuous
on $\R$
and such that $f(x)>0, g(x)>0$ for all $x>0$. Let
$(u,v)  \in C^1(\Om)\times C^1(\Om)$ be a weak solution of
\begin{equation}
\label{1.2000}\left\{\begin{array}{llll}
-\plun u=f(v)&\textrm{in }\Omega,& u>0\quad\textrm{in
}\Omega,&u=0\quad\textrm{on }\partial\Omega,\vspace*{0.2cm}\\
-\pldeux v=g(u)&\textrm{in }\Omega,& v>0\quad\textrm{in
}\Omega,&v=0\quad\textrm{on
}\partial\Omega.
\end{array}\right.
\end{equation}
The main goal of this paper is to use
the moving hyperplanes method in view of extending to a system like
(\ref{1.2000})
the monotonicity and symmetry results of
Damascelli and Pacella contained in their very nice recent article
\cite{Da}.\\

We will consider separately the cases $p_1\in(1,\infty),
 p_2=2$ (similarly
$p_2\in(1,\infty),p_1=2$) and $p_1,p_2\in(1,2)$. The first case will be
treated in a quite classical way, by using partly some comparison
principles on small
domains but also the Hopf Lemma and the strong maximum principle for the
p-Laplacian
(cf. \cite{vasq}). On the other hand
for the second case, we will establish some
monotonicity results which will be variants of some earlier
theorems
of Damascelli and Pacella in \cite{DA,Da} by using the same ideas
as in \cite{Da}, but adapted in the case of a system.\\

Before stating the monotonicity results, we
first introduce some notations used in \cite{DA,Da}.
 For any direction $\nu\in\R^N$, $|\nu|=1$, we define
$$
a(\nu):=\inf_{x\in\Omega}x.\nu,
$$
and for all $\lambda\ge a(\nu)$,
$$
\begin{array}{rcl}
\Omega_\lambda^\nu&:=&\{x\in\Omega\,|\,x.\nu<\lambda\}
(\ne\emptyset\textrm{ for }\lambda>a(\nu),\,\lambda-a(\nu)\textrm{
small}),\vspace*{0.2cm}\\
T_\lambda^\nu&:=&\{x\in\Omega\,|\,x.\nu=\lambda\}.\\
\end{array}
$$
Let us denote by $R_\lambda^\nu$ the reflection with respect to the
hyperplane $T_\lambda^\nu$ and by
$$
\begin{array}{rcl}
x_\lambda^\nu&:=&R_\lambda^\nu(x) \,\,\forall x\in\R^N,\vspace*{0.2cm}\\
(\Omega_\lambda^\nu)'&:=&R_\lambda^\nu(\Omega_\lambda^\nu),\vspace*{0.2cm}\\
\Lambda_1(\nu)&:=&\{\mu>a(\nu)\,|\,\forall\lambda\in(a(\nu),\mu),\textrm{
we do  have (\ref{situation1}) and (\ref{ii})}\},\vspace*{0.2cm}\\
\lambda_1(\nu)&:=&\sup\Lambda_1(\nu),
\end{array}
$$
where (\ref{situation1}), (\ref{ii}) are defined as follows:
\begin{equation} \label{situation1}
(\Omega_\lambda^\nu)'\textrm{ is not internally tangent to }\partial\Omega
\textrm{ at some point }p\notin T_\lambda^\nu,
\end{equation}
\begin{equation}\label{ii}
\nu(x).\nu\ne0\quad\textrm{for all }x\in\partial\Omega\cap T^\nu_\lambda,
\end{equation}
where $\nu(x)$ denotes the inward unit normal to
$\partial\Omega$ at $x$.
Notice that  since for $\lambda>a(\nu)$ and if $\lambda$ is
close to $a(\nu)$, (\ref{situation1}) and (\ref{ii}) are satisfied and
$\Omega$ is
bounded, it follows that \begin{equation}\label{53}
\Lambda_1(\nu)\ne\emptyset\quad\textrm{and}\quad\lambda_1(\nu)<\infty.
\end{equation}
Observe also that for all $\lambda>a(\nu)$, for all $c\in T_\lambda^
{\nu(x)}\cap
\Omega$ we have
\begin{equation}\label{54}
dist(c,\partial\Omega)\le\lambda-a(\nu).
\end{equation}

The monotonicity results are the following:
\begin{thm}\label{5.3}
Let $\Omega\subset\R^N$ be a bounded domain satisfying the
interior sphere condition and
let $f,g:\R\to\R^+$ be  nondecreasing on $\R^+$ and locally Lipschitz
 continuous on $\R$. Let $(u,v)  \in C^1_0(\Om)\times C^1_0(\Om)$ be a weak
solution of
$$\left\{\begin{array}{lllll}
-\plun u&=&f(v)&\textrm{in }\Omega,& u>0\quad\textrm{in
}\Omega,\vspace*{0.2cm}\\
-\Delta v&=&g(u)&\textrm{in }\Omega,& v>0\quad\textrm{in }\Omega,
\end{array}\right.
$$
where $1<p_1$. Then, for any direction $\nu\in\R^N$ and for any $\lambda$
in the interval $(a(\nu),\lambda_1(\nu)]$, we have
$$
u(x)\le u(x_\lambda^\nu)\quad\textrm{and}\quad v(x)\le v(x_\lambda^\nu)
\quad\forall  x\in\Omega_\lambda^\nu.
$$
Moreover
\begin{equation}\label{mysteres}
 \frac{\partial v}{\partial\nu}>0\quad
 \textrm{in }\Omega_\lambda^\nu\quad\forall
 \lambda<\lambda_1(\nu).
\end{equation}
\end{thm}
The following result is
the analogue of Theorem 1.1 from \cite{Da} for a system with increasing
right-hand sides.
\begin{thm}\label{mono}
Let $\Omega\subset\R^N$ be a bounded  domain with $C^1$ boundary and
let
$f,g:\R\to\R^+$ be  strictly increasing on $\R^+$, locally Lipschitz continuous
on $\R$
and such that $f(x)>0, g(x)>0$ for all $x>0$. Let $(u,v)  \in
C^1_0(\Om)\times C^1_0(\Om)$ be a weak solution of
\begin{equation}
\label{1.27}\left\{\begin{array}{lll}
-\plun u=f(v)&\textrm{in }\Omega,& u>0\quad\textrm{in }\Omega,\vspace*{0.2cm}\\
-\pldeux v=g(u)&\textrm{in }\Omega,& v>0\quad\textrm{in }\Omega,
\end{array}\right.
\end{equation}
where $p_1,p_2\in(1,2)$.
Then we have
\begin{equation}\label{rcx}
u(x)\le u(x_\lambda^\nu)\quad\textrm{and}\quad v(x)\le v(x_\lambda^\nu)
\quad\forall  x\in\Omega_\lambda^\nu,\;\forall\nu\in\R^N,\;
\forall\lambda\in(a(\nu),\lambda_1(\nu)] .
\end{equation}
\end{thm}
In Theorem \ref{mono}, the restriction $p_1,p_2\in(1,2)$ is due
to the fact that if both $p_1,p_2$
are different from 2, we must use comparison principles, and these
are less powerful if $p_1$ or $p_2$ is greater than 2. On the
other hand, if $p_1$ (or $p_2$) is equal to 2, then, as already mentioned,
we may partly
use strong maximum principles, and this finally allows $p_2$ to take
values greater or smaller than 2.
Note that this restriction is also present in the monotonicity result of
\cite{Da} in the case of a single equation.
We emphasize that in
Theorem \ref{5.3}, this condition is not needed if $p_1$ or $p_2$
is equal to two.\\
\remark \rm In \cite{Da}, Damascelli and Pacella state Theorem 1.1 under
the hypothesis
that $\Omega $ is smooth. This condition is due to the fact that
they use a sophisticated method consisting of moving
hyperplanes perpendicularly to directions $\nu$ in a neighborhood
of a fixed direction $\nu_0$. To be efficient, this method require
the continuity of $a(\nu)$ and the lower semicontinuity of $\lambda_1(\nu)$
with respect to $\nu$, and to insure this continuity
(and only for that reason),
 they assume $\Omega$
to be smooth. It appears (see \cite{Az2}) that this continuity is
guaranteed for a domain $\Omega$ of class $C^1$. To prove Theorem
\ref{mono}, we use the new technique of Damascelli and Pacella,
and so we require $\Omega$ to be $C^1$. Observe that this condition does not
appear if $p_2=2,p_1>1$. Indeed, in this case, we can use the
classical moving plane procedure consisting in moving planes
perpendicularly to a fixed direction $\nu_0$.\\

We
obtain as a consequence of Theorems \ref{mono} and \ref{5.3} the following
symmetry result:
\begin{thm}\label{corol}Let $\nu\in\R^N$ and $\Omega\subset\in\R^N$ ($N\ge2$)
be a
domain with $C^1$ boundary symmetric
with respect to the hyperplane $T_0^\nu=\{x\in\R^N\,|\,x.\nu=0\}$
and $\lambda_1(\nu)=\lambda_1(-\nu)=0$. Assume that one of the following
conditions holds:
\begin{enumerate}
\item  $p_1,p_2\in(1,2)$ and  $f,g:\R\to\R^+$ are  strictly
increasing functions on $\R^+$ such that $f(x)>0, g(x)>0$ for all $x>0$,
\item $p_1\in(1,\infty),p_2=2$ and $f,g:\R\to\R^+$ are nondecreasing
on $\R^+$.
\end{enumerate}
Moreover suppose that $f$ and $g$ are locally Lipschitz continuous
on $\R$. Then, if $(u,v)  \in C^1_0(\Om)\times C^1_0(\Om)$ is
a weak solution of (\ref{1.27})
it follows that $u$
 and $v$ are symmetric and decreasing.  In
 particular, if $\Omega$ is the ball $B_R(0)$ in $\R^N$ with center at the
origin
and radius $R$,
then $u,v$ are radially symmetric. Moreover if $f(x)>0, g(x)>0$ for all $x>0$,
then $u'(r),v'(r)<0$ for $r\in(0,R)$, $r=|x|$.
\end{thm}

Theorems \ref{mono} and \ref{5.3} have also a relatively big impact in the
study of p-Laplace systems since they are used in
\cite{Az1} to prove by blow-up some existence results and a-priori
estimates for positive solutions of the system
\begin{equation}\label{azerty}
\left\{\begin{array}{llll}
-\plun u=f(|v|)&\textrm{ in }\Omega,&u=0&\textrm{
on }\partial\Omega,\vspace*{0.3cm}\\
-\pldeux v=g(|u|)&\textrm{ in }\Omega,&v=0&\textrm{
on }\partial\Omega,
\end{array}\right.
\end{equation}
where $1<p_1,p_2<N$, $\Omega$ is convex,
$f,g:\R\to\R^+$ are nondecreasing locally Lipschitz
continuous on $(0,+\infty)$, continuous on $[0,+\infty)$ and
satisfy
\begin{equation}\label{cdt}
C_1|s|^{q_1}\le f(s)\le C_2|s|^{q_1},\quad D_1|s|^{q_2}\le g(s)\le
D_2|s|^{q_2}\quad\forall s\in\R^+
\end{equation}
for some positive constants $C_1,C_2,D_1,D_2$ and
$q_1q_2>(p_1-1)(p_2-1)$.\\

This paper is organized as follows. In section \ref{prelim}, we
recall some well known results concerning the p-Laplacian
operator.
In section \ref{341}, we prove some weak comparison principles on small domains
which
are some adaptations of
Theorem 1.2 from \cite{DA} to systems. In section \ref{342},
we use these principles to prove the
monotonicity results and finally, we prove as a corollary Theorem
\ref{corol}.\\
\mbox{}\\
\textbf{Acknowledgments}\\
I would like to acknowledge particularly Professors Philippe
Cl\'ement and Enzo Mitidieri for many useful comments.
\section{Preliminaries}\label{prelim}\setcounter{equation}{0}
Suppose that $f,g$ are given positive continuous functions as in the
introduction.
\begin{definition}
\rm Let $t\ge0$. A function $(u,v)\in C^1_0(\Om)\times C^1_0(\Om)$
is said a weak solution of (\ref{1.2000}) if for any function
$\phi\in C^\infty_c(\Omega)$
we have
\begin{equation}
\label{faible}
\left\{\begin{array}{l}
\displaystyle\io|\nabla u|^{p_1-2}\nabla u.\nabla\phi\, dx=\io f(v)\phi\,
dx,\\
\displaystyle\io|\nabla u|^{p_1-2}\nabla u.\nabla\phi\, dx=\io g(u)\phi\,
dx.
\end{array}\right.
\end{equation}
\end{definition}
We are interested in monotonicity results for weak solutions of
(\ref{1.2000}).\\

By the maximum principle and Hopf's lemma
of \cite{vasq} for the
p-Laplacian, any weak solution $(u,v)$ of (\ref{1.2000})
satisfies
$$
u>0\quad\textrm{in }\Omega,\quad\frac{\partial
u}{\partial\nu}<0\quad\textrm{on }\partial\Omega,$$
$$
v>0\quad\textrm{in }\Omega,\quad\frac{\partial
v}{\partial\nu}<0\quad\textrm{on }\partial\Omega
$$
where $\nu$ denotes the outward unit normal to $\partial\Omega$.\\

In the present section we recall some well known properties of the
operator $-\pl$. The following result is due to Damascelli (\cite{DA}).

\begin{lemma}[Weak comparison principle]\label{compa}
Let $p>1$. If
$u,v\in W^{1,\infty}(\Omega)$ are such that
\begin{equation}
\label{compar1}
\io|\nabla u|^{p-2}\nabla u.\nabla\phi\,dx\le\io|\nabla v|^{p-2}\nabla
v.\nabla\phi
\,dx\qquad\forall\phi\in C^\infty_c(\Omega),\phi\ge0
\end{equation}
and $u\le v$ on
$\partial\Omega$, then
$u\le v$ on $\Omega$.
\end{lemma}
Next we state a strong comparison principle due to
Damascelli in \cite{DA} (Theorem 1.4).
\begin{lemma}[Strong comparison principle]
\label{thm1.4}
Let $\Omega\subset\R^N$ be a bounded domain and $p>1$.
Let $u,v\in C^1(\Omega)$ satisfy
$$\left\{\begin{array}{l}
\displaystyle\io|\nabla u|^{p-2}\nabla u.\nabla\phi\,dx\le\io|\nabla
v|^{p-2}\nabla v.\nabla\phi
\,dx\qquad\forall\phi\in C^\infty_c(\Omega),\phi\ge0,\\
u\le v\quad in\quad\Omega
\end{array}\right.$$ and define $Z:=\{x\in\Omega\,|\,|\nabla u(x)|+|\nabla
v(x)|=0\}$
if $p\ne2$, $Z:=\emptyset$ if $p=2$.

If $x_0\in\Omega\setminus Z$ and $u(x_0)=v(x_0)$, then $u=v$ in
the connected component of $\Omega\setminus Z$ containing $x_0$.
\end{lemma}
Finally we recall a lemma proved by Simon in \cite{simon} and
Damascelli in \cite{DA} which will be used later.
\begin{lemma}
\label{simon}
Let $p>1$ and $N\in\N_0$. There exist some positive constants $c_1,c_2$
depending on p and $N$ such that
for all $\eta,\eta'\in\R^N$ with $|\eta|+|\eta'|>0$
\begin{equation}
\label{simon1}
||\eta|^{p-2}\eta-|\eta'|^{p-2}\eta'|\le
c_1(|\eta|+|\eta'|)^{p-2}|\eta-\eta'|
\end{equation}
\begin{equation}
\label{simon2}
(|\eta|^{p-2}\eta-|\eta'|^{p-2}\eta').(\eta-\eta')\ge
c_2(|\eta|+|\eta'|)^{p-2}|\eta-\eta'|^2.
\end{equation}
\end{lemma}
\section{Weak comparison principles}\label{341}

Let $\Omega$ be a bounded domain in $\R^N$  and let $(u,v), (\bar u,\bar v)
\in \espaces$
be solutions of
\begin{equation}\label{etoiles}
\left\{\begin{array}{ll}
-\plun u=f(v)&\textrm{on}\;\Omega,\quad u\ge0\vspace*{0.2cm}\\
-\pldeux v=g(u)&\textrm{on} \;\Omega,\quad v\ge0
\end{array}\right.
\end{equation}
\begin{equation}\label{stars}
\left\{\begin{array}{ll}
-\plun \bar u=f(\bar v)&\textrm{on}\;\Omega,\quad \bar u\ge0\vspace*{0.2cm}\\
-\pldeux \bar v=g(\bar u)&\textrm{on} \;\Omega,\quad \bar v\ge0
\end{array}\right.
\end{equation}
where $t>0$ is a real parameter and $f,g:\R\to\R^+$ are  locally Lipschitz
continuous
on $\R$  and
nondecreasing on $\R^+$. As mentioned above, our first aim is to
prove some comparison principles for solutions of
(\ref{etoiles}), (\ref{stars}).
We begin with a result in case $p_1,p_2\in(1,2)$
that
will be an extension of Theorem
2.2 of \cite{Da} to  systems with nondecreasing right
hand side.\\

 For any set $A\subseteq\Omega$, we define $M_A=M_A(u,\bar
u):=\sup_A(|\nabla u|+|\nabla \bar u|)$ and $\tilde M_A=\tilde M_A(v,\bar
v):=\sup_A(|\nabla v|+|\nabla \bar v|)$. We shall denote the measure
of a measurable set $B$ by $|B|$.
\begin{thm}\label{thm2.2}
Let $\Omega$ be a bounded domain contained in $\R^N$ and suppose
that $1<p_1,p_2<2$. Let $(u,v),(\bar u,\bar v)\in C^1(\Om)\times C^1(\Om)$
be two solutions of (\ref{etoiles}), (\ref{stars}). Suppose that
$$f,g:\R\to\R^+$$
are  locally Lipschitz continuous on $\R$  and
nondecreasing on $\R^+$. Then there exist $\alpha,M>0$ depending
on $N,p_1,p_2$, $f,g$,
$|\Omega|$, $M_\Omega$, $\tilde M_\Omega$ and the $L^\infty$
norms of $u,\bar u,v,\bar v$ such that if $\Omega'\subset\Omega$
is an open set and if there exists measurable sets $A_i,\tilde A_i (i=1,2,3)$
such that
$$\Omega'=A_1\cup A_2\cup A_3=\tilde A_1\cup \tilde A_2\cup \tilde
A_3$$ and
$$\left\{\begin{array}{l}
A_1\cup  A_2, \tilde A_1\cup \tilde A_2\quad\textrm{are open}\vspace*{0.2cm}\\
A_i\cap A_j=\tilde A_i\cap \tilde A_j=\emptyset\quad\textrm{for
all } i\ne j
\end{array}\right.$$
and
$$\left\{\begin{array}{l}
\max\{|A_1|,|\tilde A_1|\}<\alpha,\vspace*{0.2cm}\\
\max\{M_{A_2},\tilde M_{\tilde
A_2}\}<M,\vspace*{0.2cm}\\
u< \bar{u}\quad\textrm{on }A_3,\vspace*{0.2cm}\\
v<\bar v\quad\textrm{on }\tilde A_3,
\end{array}\right.$$
then we have the implication
$$
\left\{\begin{array}{l}
u\le\bar u \textrm{ on }\partial\Omega'\cup\partial(A_1\cup A_2)\\
v\le\bar v \textrm{ on }\partial\Omega'\cup\partial(\tilde A_1\cup\tilde A_2)
\end{array}\right.\Longrightarrow\left\{\begin{array}{l}
u\le\bar u\\v\le\bar v
\end{array}\right.\textrm{ on }\Omega'.$$
\end{thm}
\proof By multiplying
the first equations of (\ref{etoiles}), (\ref{stars}) by $(u-\bar u)^+\in
W^{1,p_1}_0(\Omega')$
(cf.. Theorem IX.17. and remark 20, p. 171-172 in \cite{Brezis}) and the
second equations
by $(v-\bar v)^+\in W^{1,p_2}_0(\Omega')$, and subtracting
the resulting identities, we get
\begin{equation}\label{40}
\int_{\Omega'\cap[u\ge\bar u]}(|\nabla u|^{p_1-2}\nabla u-|\nabla \bar
u|^{p_1-2}\nabla \bar
u).\nabla (u-\bar u)\,dx=\int_{\Omega'\cap[u\ge\bar
u]}(f(v)-f(\bar v))(u-\bar u)\,dx,
\end{equation}
\begin{equation}
\label{41}
\int_{\Omega'\cap[v\ge\bar v]}(|\nabla v|^{p_2-2}\nabla v-|\nabla \bar
v|^{p_2-2}\nabla \bar
v).\nabla (v-\bar v)\,dx=\int_{\Omega'\cap[v\ge\bar
v]}(g(u)-g(\bar u))(v-\bar v)\,dx.
\end{equation}
By Lemma \ref{simon}, the left-hand sides of (\ref{40}), (\ref{41}) are
respectively
greater or equal to
\begin{equation}
\label{415}
c_2M_\Omega^{p_1-2}\int_{A_1\cap[u\ge\bar u]}|\nabla(u-\bar
u)|^2\,dx+c_2M_{A_2}^{p_1-2}\int_{A_2\cap[u\ge\bar u]}|\nabla(u-\bar
u)|^2\,dx
\end{equation}
and
$$c_2\tilde M_\Omega^{p_2-2}\int_{\tilde A_1\cap[v\ge\bar v]}|\nabla(v-\bar
v)|^2\,dx+c_2\tilde M_{\tilde A_2}^{p_2-2}\int_{\tilde A_2\cap[v\ge\bar v]}
|\nabla(v-\bar
v)|^2\,dx$$
where $c_2$ is a positive constant depending on $p_1,p_2$ and $N$.
Since $f$ and $g$ are nondecreasing, the right-hand side of (\ref{40})
can be further majorized with
$$\int_{\Omega'\cap[u\ge\bar
u]\cap[v\ge\bar v]}(f(v)-f(\bar v))(u-\bar u)\,dx,$$
and by the local Lipschitz property of $f,g$, this latter quantity
is smaller or equal to
\begin{eqnarray*}
\Lambda\int_{\Omega'\cap[u\ge\bar
u]\cap[v\ge\bar v]}(v-\bar v)(u-\bar u)&\le&\Lambda\|(v-\bar v)^+
\|_{L^2((\tilde A_1\cup \tilde A_2)\cap[u\ge \bar u])}\|(u-\bar u)^+
\|_{L^2((A_1\cup A_2)\cap[v\ge \bar v])}\\
{}&\le &\Lambda\|(v-\bar v)^+
\|_{L^2((\tilde A_1\cup \tilde A_2))}\|(u-\bar u)^+
\|_{L^2((A_1\cup A_2))},
\end{eqnarray*}
for some constant $\Lambda>0$ depending on $f$ and the $L^\infty$ norms of
$v,\bar v$ (cf.. remark 2.1 in \cite{Da}). Using a version of Poincar\'e's
inequality
(see Lemma 2.2 of \cite{DA}),
this last term is smaller than
$$\Lambda w_N^{-2/N}|\Omega'|^{1/N}\left\{|\tilde A_1|^{1/2N}\|\nabla
(v-\bar v)
\|_{L^2(\tilde A_1\cap[v\ge\bar v])}+|\Omega|^{1/2N}\|\nabla (v-\bar v)
\|_{L^2(\tilde A_2\cap[v\ge\bar v])}\right\}$$
$$
\times\left\{|A_1|^{1/2N}\|\nabla (u-\bar u)
\|_{L^2( A_1\cap[u\ge\bar u])}+|\Omega|^{1/2N}\|\nabla (u-\bar u)
\|_{L^2( A_2\cap[u\ge\bar u])}\right\},$$
where $\omega_N$ denotes the Lebesgue measure of the unit ball in $\R^N$.
The same reasoning can be made with (\ref{41}) (with another constant
$\Lambda'$).
Adding both inequalities, we obtain
\begin{eqnarray}
{}&{}&\nonumber c_2M_\Omega^{p_1-2}\|\nabla (u-\bar u)
\|^2_{L^2( A_1\cap[u\ge\bar u])}+c_2M_{A_2}^{p_1-2}\|\nabla (u-\bar u)
\|^2_{L^2( A_2\cap[u\ge\bar u])}
\\ {}&{}&\nonumber +c_2\tilde M_\Omega^{p_2-2}\|\nabla (v-\bar v)
\|^2_{L^2(\tilde A_1\cap[v\ge\bar v])}
+c_2\tilde M_{\tilde A_2}^{p_2-2}\|\nabla (v-\bar v)
\|^2_{L^2(\tilde A_2\cap[v\ge\bar v])}
\\ {}&{}&\nonumber
\le 2\max\{\Lambda,\Lambda'\}w_N^{\frac{-2}{N}}|\Omega'|^\frac{1}{N}\left\{
|A_1|^\frac{1}{2N}|\tilde A_1|^\frac{1}{2N}\|\nabla (v-\bar v)
\|_{L^2(\tilde A_1\cap[v\ge\bar v])}\|\nabla (u-\bar u)
\|_{L^2( A_1\cap[u\ge\bar u])}\right.\\ {}&{}&\nonumber
+|\Omega|^\frac{1}{N}\|\nabla (v-\bar v)
\|_{L^2(\tilde A_2\cap[v\ge\bar v])}\|\nabla (u-\bar u)
\|_{L^2( A_2\cap[u\ge\bar u])}\\ {}&{}&\nonumber
+|A_1|^\frac{1}{2N}|\Omega|^\frac{1}{2N}\|\nabla (v-\bar v)
\|_{L^2(\tilde A_2\cap[v\ge\bar v])}\|\nabla (u-\bar u)
\|_{L^2( A_1\cap[u\ge\bar u])}\\ {}&{}&\label{42} \left.
+|\Omega|^\frac{1}{2N}|\tilde A_1|^\frac{1}{2N}\|\nabla (v-\bar v)
\|_{L^2(\tilde A_1\cap[v\ge\bar v])}\|\nabla (u-\bar u)
\|_{L^2( A_2\cap[u\ge\bar u])}\right\}.
\end{eqnarray}
By Young inequality, the right-hand side of
(\ref{42}) is smaller or equal to
\begin{eqnarray}{}&{}&\nonumber
2\max\{\Lambda,\Lambda'\}w_N^{\frac{-2}{N}}|\Omega'|^\frac{1}{N}\left\{|\tilde
A_1|^\frac{1}{N}\|\nabla (v-\bar v)
\|^2_{L^2(\tilde A_1\cap[v\ge\bar v])}+|A_1|^\frac{1}{N}\|\nabla (u-\bar u)
\|^2_{L^2( A_1\cap[u\ge\bar u])}\right.\\ \nonumber
{}&{}&\left.+|\Omega|^\frac{1}{N}\|\nabla (u-\bar u)
\|^2_{L^2( A_2\cap[u\ge\bar u])}+|\Omega|^\frac{1}{N}
\|\nabla (v-\bar v)
\|^2_{L^2(\tilde A_2\cap[v\ge\bar v])}\right\}.
\end{eqnarray}
From this we infer that if $|A_1|,|\tilde A_1|,M_{A_2}$ and $\tilde
M_{\tilde A_2}$
are small enough, then
$$\begin{array}{ll}
{}&\|\nabla (v-\bar v)
\|_{L^2(\tilde A_1\cap[v\ge\bar v])}=\|\nabla (u-\bar u)
\|_{L^2( A_1\cap[u\ge\bar u])}\\=&\|\nabla (u-\bar u)
\|_{L^2( A_2\cap[u\ge\bar u])}=\|\nabla (v-\bar v)
\|_{L^2(\tilde A_2\cap[v\ge\bar v])}=0,\end{array}$$
so that by Poincar\'e's inequality,  $(u-\bar u)^+=(v-\bar v)^+=0$
in respectively $A_1\cup A_2$ and $\tilde A_1\cup\tilde A_2$, and finally on
$\Omega'$ by definition of $A_3,\tilde A_3$.
\cqfd
We now give a weak comparison principle in the case $p_1\in(1,+\infty),p_2=2$
(or $p_1=2,p_2\in(1,+\infty)$). Let $u,\bar u,v,\bar v\in
C^1(\Om)$ and $A\subset\Omega$. Set
$$m_A:=\inf_A(|\nabla u|+|\nabla \bar u|),\quad\tilde m_A:=\inf_A(|\nabla v|+
|\nabla \bar v|).$$
\begin{thm}
\label{wcp}
Let $m>0$, $\Omega\subset\R^N$ be a bounded domain and $(u,v),(\bar u,\bar
v)\in\espaces$
be two solutions of (\ref{etoiles}), (\ref{stars}) where
$f,g:\R\to\R^+$ are  nondecreasing on $\R^+$ and
locally Lipschitz continuous on $\R$.
\begin{enumerate}
\item[(i)] If $p_1\in(1,2),p_2=2$, then there exists $\delta>0$
depending on $N,p_1, M_\Omega,f,g$ and the $L^\infty$-norms
of $u,\bar u,v,\bar v$ such that for any open subset $\Omega'\subseteq\Omega$
with $|\Omega'|<\delta$, $u\le\bar u$ and $v\le\bar v$ on $\partial\Omega'$
implies $u\le\bar u$ and $v\le\bar v$ on $\Omega'$.
\item[(ii)] If $p_1>2,p_2=2$,
then
there exists $\delta>0$   depending on $N,p_1, m,f,g$ and
the $L^\infty$ norms of $u,\bar u,v,\bar v$
such that if such that if $m_\Omega\ge m$, if $\Omega'\subseteq\Omega$ is
an open subset with
$|\Omega'|<\delta$, then $u\le\bar u$ and $v\le\bar v$ on $\partial\Omega'$
implies $u\le\bar u$ and $v\le\bar v$ on $\Omega'$.
\item [(iii)] If $p_1=p_2=2$, then there exists $\delta>0$
depending on $N,f,g$ and the $L^\infty$-norms
of $u,\bar u,v,\bar v$ such that for any open subset $\Omega'\subseteq\Omega$
with $|\Omega'|<\delta$, $u\le\bar u$ and $v\le\bar v$ on $\partial\Omega'$
implies $u\le\bar u$ and $v\le\bar v$ on $\Omega'$.
\end{enumerate}
\end{thm}
\proof
Let us prove (i). As in the proof of Theorem \ref{thm2.2}, we
first
multiply the first equations of (\ref{etoiles}), (\ref{stars}) by
$(u-\bar u)^+\in W^{1,p_1}_0(\Omega')$ and
 the second equations by $(v-\bar v)^+\in W^{1,2}_0(\Omega')$
and we subtract them. In this way we obtain (\ref{40}), while (\ref{41}) is
replaced by
\begin{equation}\label{400}
  \int_{\Omega'\cap[v\ge \bar v]}|\nabla(v-\bar v)|^2\,dx=
  \int_{\Omega'\cap[v\ge \bar
  v]}(g(u)-g(\bar u))(v-\bar
  v)\,dx.
\end{equation}
By Lemma \ref{simon}, the left-hand side of (\ref{40}) can be
estimated with
$$
c_2M_\Omega^{p_1-2}\int_{\Omega'\cap[u\ge \bar u]}|\nabla(u-\bar u)|^2\,dx.$$
We then treat the right-hand sides of (\ref{40}), (\ref{400}) as in
the proof of Theorem \ref{thm2.2}. For this purpose we use Lemma 2.2 from
\cite{DA} with $A_1=\Omega',A_2=\emptyset$ and
we add the obtained inequalities
 to get
$$\begin{array}{l}
c_2M_\Omega^{p_1-2}\|\nabla(u-\bar u)^+\|^2_{L^2(\Omega')}+\|\nabla(v-\bar v)^+
\|^2_{L^2(\Omega')}\\\le4\max\{\Lambda,\Lambda'\}\omega_N^{-\frac{2}{N}}|\Omega'
|^\frac{2}{N}(\|\nabla(u-\bar u)^+\|^2_{L^2(\Omega')}+\|\nabla(v-\bar v)^+
\|^2_{L^2(\Omega')}).
\end{array}$$
So, if $|\Omega'|$ is sufficiently small, then necessarily
$\|\nabla(u-\bar u)^+\|_{L^2(\Omega')}=\|\nabla(v-\bar v)^+
\|_{L^2(\Omega')}=0$. By Poincar\'e's inequality this implies $(u-\bar
u)^+=(v-\bar v)^+=0$
in $\Omega'$, i.e. $u\le\bar u$ and $v\le \bar v$ in $\Omega'$.\\

Let us now prove (ii). We again obtain (\ref{40}), (\ref{400}) and by
Lemma \ref{simon}, the right-hand side of (\ref{40}) is greater or
equal to
$$
c_2 m_{\Omega}^{p_1-2}\int_{\Omega'\cap[u\ge \bar u]}
|\nabla(u-\bar u)|^2\,dx.$$
Writing as above the same estimates of the left-hand sides of (\ref{40}),
(\ref{400}), we get:
$$\begin{array}{l}
c_2m_{\Omega}^{p_1-2}\|\nabla(u-\bar
u)^+\|^2_{L^2(\Omega')}+\|\nabla(v-\bar v)^+
\|^2_{L^2(\Omega')}\\\le4\max\{\Lambda,\Lambda'\}\omega_N^{-\frac{2}{N}}|\Omega'
|^\frac{2}{N}(\|\nabla(u-\bar u)^+\|^2_{L^2(\Omega')}+\|\nabla(v-\bar v)^+
\|^2_{L^2(\Omega')}),
\end{array}$$
and we conclude as in the proof of (i).\\

Finally, we prove (iii).
We get (\ref{400}) while (\ref{40}) is replaced by
\begin{equation}\label{4000}
  \int_{\Omega'\cap[u\ge \bar u]}|\nabla(u-\bar u)|^2\,dx=
  \int_{\Omega'\cap[u\ge \bar
  u]}(f(v)-f(\bar v))(u-\bar
  u)\,dx.
\end{equation}
We estimates as in the proof of (i) the left-hand sides of
(\ref{4000}), (\ref{400}) and we obtain
$$\begin{array}{l}
\|\nabla(u-\bar u)^+\|^2_{L^2(\Omega')}+\|\nabla(v-\bar v)^+
\|^2_{L^2(\Omega')}\\\le4\max\{\Lambda,\Lambda'\}\omega_N^{-\frac{2}{N}}|\Omega'
|^\frac{2}{N}(\|\nabla(u-\bar u)^+\|^2_{L^2(\Omega')}+\|\nabla(v-\bar v)^+
\|^2_{L^2(\Omega')}).
\end{array}$$
We conclude again as in the proof of (i).\cqfd


\subsection{Proof of the monotonicity results}\label{342}

 Let us introduce some more notations used in \cite{Da}.
For any direction $\nu\in\R^N$, $|\nu|=1$, we define
$$\Lambda_2(\nu)=\{\lambda>a(\nu)\,|\,(\Omega_\mu^\nu)'\subset\Omega\textrm{ for
any }
\mu\in(a(\nu),\lambda]\},$$
and, if $\Lambda_2(\nu)\ne\emptyset$,
$$\lambda_2(\nu)=\sup\Lambda_2(\nu).$$ If $a(\nu)<\lambda\le
\lambda_2(\nu)$, $x\in\Omega_\lambda^\nu$, $u,v\in C^1(\Om)$, we set
$$u_\lambda^\nu(x)=u(x_\lambda^\nu),\quad v_\lambda^\nu(x)=v(x_\lambda^\nu),$$
$$Z_\lambda^\nu=Z_\lambda^\nu(u)=\{x\in\Omega_\lambda^\nu\,|
\,\nabla u(x)=\nabla u_\lambda^\nu(x)=0\},$$
$$\tilde Z_\lambda^\nu=\tilde Z_\lambda^\nu(v)=\{x\in\Omega_\lambda^\nu\,|
\,\nabla v(x)=\nabla v_\lambda^\nu(x)=0\}$$ and
$$Z=Z(u)=\{x\in\Omega\,|\,\nabla u(x)=0\},$$
$$\tilde Z=\tilde Z(v)=\{x\in\Omega\,|\,\nabla v(x)=0\}.$$
We also define
$$\Lambda_0(\nu)=\{\lambda\in (a(\nu),\lambda_2(\nu)]\,|\,u\le u_\mu^\nu
\textrm{ and }v\le v_\mu^\nu\textrm{ in }\Omega_\mu^\nu
\textrm{ for any }\mu\in(a(\nu),\lambda]\}$$
and if $\Lambda_0(\nu)\ne\emptyset$, we set
$$\lambda_0(\nu)=\sup\Lambda_0(\nu). $$
As remarked in \cite{Da}, we obviously have
$\lambda_0(\nu)\le\lambda_1(\nu)\le\lambda_2(\nu)$.\\

We begin now to prove Theorem \ref{5.3}. In the proof we shall use the weak
comparison principles
stated in
Theorem \ref{wcp} for the beginning of the moving plane procedure,
but afterwards it becomes quite classical in the sense that it uses
maximum principles and Hopf's lemma for the usual Laplacian. That's why the
result is true for all
$p_2\in (0,+\infty)$, in opposition with Theorem \ref{mono}.\\
\mbox{}\\
\textbf{Proof of Theorem \ref{5.3}:} Let us fix a direction $\nu$. If
$\lambda\le\lambda_1(\nu)$, we have $u\le u_\lambda^\nu$,
$v\le v_\lambda^\nu$ on $\partial\Omega_\lambda^\nu$ since $u,v>0$ on $\Omega$,
$u=v=0$ on
$\partial\Omega$. If $p_1\le 2$,
there exists $\delta>0$ such that Theorem \ref{wcp} (i) or (iii) is
applicable to
the pairs $(u,v)$, $(\bar u,\bar v)=(u_\lambda^\nu,v_\lambda^\nu)$ and
$\Omega'=\Omega_\lambda^\nu$ for all $\lambda\in(a(\nu),\lambda_2(\nu))$.
Since for $\lambda>a(\nu)$,
$\lambda-a(\nu)$ small enough, we have
$\left|\Omega_\lambda^\nu\right|<\delta$,
we get
$u\le u_\lambda^\nu$ and $v\le v_\lambda^\nu$
on $\Omega_\lambda^\nu$ for these values of $\lambda$.\\
If $p_1>2$, then, by the Hopf's lemma (see \cite{vasq}), there exists
$\bar\lambda>a(\nu)$ and $m>0$ such that  $m_{\Omega_\lambda^\nu}\ge m$
for all $\lambda\in(a(\nu),\bar\lambda)$. Moreover, as above, for
$\lambda>a(\nu)$,
$\lambda-a(\nu)$ small enough, it follows that
$\left|\Omega_\lambda^\nu\right|$ is small.
So for these values of $\lambda$, we can  apply Theorem \ref{wcp} (ii) with
$u,\bar u:=u_\lambda^\nu,v,\bar v:=v_\lambda^\nu$,
$\Omega=\Omega'=\Omega_\lambda^\nu$, to get
$u\le u_\lambda^\nu$ and $v\le v_\lambda^\nu$
on $\Omega_\lambda^\nu$ . This proves that $\Lambda_0(\nu)\ne\emptyset$.\\

Suppose by contradiction that $\lambda_0(\nu)<\lambda_1(\nu)$.
By the continuity of $u,v$ it follows that $u\le u_{\lambda_0(\nu)}^\nu$ and
$v\le v_{\lambda_0(\nu)}^\nu$ on $\Omega_{\lambda_0(\nu)}^\nu$. Thus, by
Lemma \ref{thm1.4}, we have either $v< v_{\lambda_0(\nu)}^\nu$ or
$v= v_{\lambda_0(\nu)}^\nu$ on $\Omega_{\lambda_0(\nu)}^\nu$. Since
 (\ref{situation1}) holds with $\lambda=\lambda_0(\nu)$, we have
$0=v<v_{\lambda_0(\nu)}^\nu$ on $\partial
\Omega_{\lambda_0(\nu)}^\nu\setminus T_{\lambda_0(\nu)}^\nu$, so that
we are in
the first case.
Using the fact that $g$ is nondecreasing, we obtain
$$\left\{\begin{array}{ll}-\Delta(v-v_{\lambda_0(\nu)}^\nu)\le0&\textrm{ on }
\Omega_{\lambda_0(\nu)}^\nu,\vspace*{0.2cm}\\
v-v_{\lambda_0(\nu)}^\nu<0 &\textrm{ on }
\Omega_{\lambda_0(\nu)}^\nu,\vspace*{0.2cm}\\
v=v_{\lambda_0(\nu)}^\nu&\textrm{ on
}T_{\lambda_0(\nu)}^\nu.\end{array}\right.$$
Thus by Hopf's lemma we have
$$
\frac{\partial(v-v_{\lambda_0(\nu)}^\nu)}{\partial\nu}>0\textrm{ on }
T_{\lambda_0(\nu)}^\nu\cap\Omega.
$$
Hence $\frac{\partial v}{\partial\nu}>0$ on
$T_{\lambda_0(\nu)}^\nu\cap\Omega$. Since
$\lambda_0(\nu)<\lambda_1(\nu)$, we have $\nu(x).\nu>0$ for all
$x\in\partial\Omega\cap\partial
\Omega_{\lambda_0(\nu)}^\nu$.  By
Hopf's lemma, for any $x\in\partial\Omega\cap\partial
\Omega_{\lambda_0(\nu)}^\nu$, it holds that
$\nabla v(x)=c(x)\nu(x)$ for some function $c(x)>0$, so that
$\frac{\partial v}{\partial\nu}>0$ on $\partial\Omega\cap\partial
\Omega_{\lambda_0(\nu)}^\nu$, and hence
\begin{equation}\label{rx}
\frac{\partial v}{\partial\nu}>0
\textrm{ on
}T_{\lambda_0(\nu)}^\nu\cap\overline{\Omega_{\lambda_0(\nu)}^\nu}.
\end{equation}
Since we suppose by contradiction that
$\lambda_0(\nu)<\lambda_1(\nu)$, by definition of $\lambda_0(\nu)$ we infer
that there exists a sequence
$(\lambda_n)\subset(\lambda_0(\nu),\lambda_1(\nu))$ such that
$\lambda_n\to\lambda_0(\nu)$
and a sequence $(x_n)\subset\Omega_{\lambda_n}^\nu$ with the
property
\begin{equation}\label{tr}
v(x_n)>v_{\lambda_n}^\nu(x_n)\quad\forall n\in\N.
\end{equation}
Indeed, if we would have $v\le v_{\lambda_n}^\nu$ in
$\Omega_{\lambda_n}^\nu$, then, due to the fact that
$f$ is nondecreasing, we would have
$u\le
u_{\lambda_n}^\nu$ in $\Omega_{\lambda_n}^\nu$. However this contradicts
$\lambda_n>\lambda_0(\nu)$. Since $(x_n)$ is bounded, there
exists $x\in\overline{\Omega_{\lambda_0(\nu)}^\nu}$ such that $x_n\to
x$. Passing to the limit as $n\to\infty$ in (\ref{tr}) we obtain
$v(x)\ge v_{\lambda_0(\nu)}^\nu(x)$. As a consequence it follows that
$x\in T_{\lambda_0(\nu)}^\nu$. By (\ref{tr}), there exists a
sequence $(y_n)\subset (x_n,(x_n)_{\lambda_n}^\nu)$ (where $(a,b)$
for $a,b\in\R^N$ denotes here the open segment of extremities $a$ and $b$)
such that $\frac{\partial v}{\partial\nu}(y_n)<0$. Clearly $y_n\to
x$ and then $\frac{\partial v}{\partial\nu}(x)\le 0$,
which contradicts (\ref{rx}).\\

To prove (\ref{mysteres}), it suffices to apply the same reasoning
as above to the function $v-v_\lambda^\nu$ on
$\Omega_\lambda^\nu$. Indeed we have
$-\Delta(v-v_\lambda^\nu)\le0$, $v< v_\lambda^\nu$
on $\Omega_\lambda^\nu$ and $v=v_\lambda^\nu$ on
$T_\lambda^\nu\cap\Omega$. So we obtain
$\frac{\partial(v-v_\lambda^\nu)}{\partial\nu}>0$ on $T_\lambda^\nu\cap\Omega$
and finally $\frac{\partial v}{\partial\nu}>0$ on
$T_\lambda^\nu\cap\Omega$. Since
$\Omega_{\lambda_0(\nu)}^\nu=
\cup_{a(\nu)<\lambda<\lambda_0(\nu)}(T_\lambda^\nu\cap\Omega)$, we
have (\ref{mysteres}). This completes the proof.
\cqfd
\mbox{}\\

Now we treat the case $p_1,p_2\in (1,2)$.
To prove Theorem \ref{mono}, we will follow the same
steps as in \cite{Da}, but adapted to our case.
For the sake of completeness and clarity, we will sometimes repeat
some arguments from \cite{Da}.\\

Since the proof is quite long, we would like to
give the main ideas beyond it.
We first prove Lemma \ref{3.1}, an extension to our system
of Theorem 3.1 from \cite{Da} (see also Theorem 1.5 from
\cite{DA}). It asserts that once we start the moving plane
procedure along a direction $\nu$, if
$\lambda_0(\nu)<\lambda_2(\nu)$, then the set $Z$ of critical points
of $u$ creates a connected component $C$ of $\Omega\setminus Z$
symmetric with respect to $T_{\lambda_0(\nu)}^\nu$ and where
$u=u_{\lambda_0(\nu)}^\nu$, and the same results holds for the
function $v$ with a component $\tilde C$ of $\Omega\setminus \tilde Z$.
Hence our goal is to
prove that such  sets $C$ or $\tilde C$ cannot
exist.

A first step in that way is Lemma \ref{prop 3.1} which implies
that if $C$ is defined as above and if $u$ is constant on a
connected subset of $\partial C$ whose projection on $T_{\lambda_0(\nu)}^\nu$
contains a relatively open nonempty subset, then such a set $C$
cannot exists, and the analogue holds for function $v$ with a
component $\tilde C$.

In fact, Lemmas \ref{3.1} and \ref{prop 3.1} imply that if $(u,v)$ is solution
of (\ref{1.27}) and if $\lambda_0(\nu)<\lambda_2(\nu)$, then there
exist a connected component $C^\nu$ of
$\Omega_{\lambda_0(\nu)}^\nu\setminus Z_{\lambda_0(\nu)}^\nu$
and a component $\tilde C^\nu$ of $\Omega_{\lambda_0(\nu)}^\nu\setminus \tilde
Z_{\lambda_0(\nu)}^\nu$ such that $u=u_{\lambda_0(\nu)}^\nu$ in $C^\nu$
and $v=v_{\lambda_0(\nu)}^\nu$ in $\tilde C^\nu$. So if we suppose
moreover that either $u$ or $v$ is constant on each connected
component of respectively $Z$ or $\tilde Z$, then in the first
case, $\partial C^\nu$ would contain a set $\Gamma$ on which $\nabla
u=0$, $u$ is constant and whose projection on $T_{\lambda_0(\nu)}^\nu$
contains an open subset of $T_{\lambda_0(\nu)}^\nu$, which will be
impossible by Lemma \ref{prop 3.1}.

At first sight, one could think that the assumption that $u$ or
$v$ is constant on respectively $Z$ or $\tilde Z$ is always
satisfied by any $C^1$ function. In fact, this is not true and the
question of finding sufficient condition on a connected set of
critical points of a $C^1$ function ensuring that this function is
constant there seems very complicated. Some counterexamples are
cited in \cite{Da} (see \cite{W} and \cite{N}), and they show that
this question is strictly related to Sard's lemma and the theory
of fractal sets.

So, to prove that either $u$ or $v$ is constant on a connected set
of respectively $\partial C$ or $\partial \tilde C$ (where $C$ and $\tilde
C$ are introduced
above) and that the projection of $C$
on the hyperplane $T_{\lambda_0(\nu)}^\nu$
contains a nonempty open set, some extra work is needed. To do
that, we use the new argument introduced by Damascelli and Pacella
in \cite{Da} consisting in moving hyperplanes orthogonal to
directions close to $\nu$ to prove that the set $C$ (or $\tilde C$) is also
symmetric with nearby hyperplanes and to show that on its
boundary, there is at least one connected piece where $u$ is
constant, $\nabla u=0$ and whose projection on $T_{\lambda_0(\nu)}^\nu$
contains a nonempty
 open set.\\

We first give the analogue for our system of Theorem 3.1 from
\cite{Da}.
\begin{lemma}
\label{3.1}Let $\Omega\subset\R^N$
be a bounded domain and
let $(u,v)\in C^1_0(\Om)\times C^1_0(\Om)$ be a weak solution of
(\ref{etoiles})
 where
$p_1,p_2\in (1,2)$, $f,g$ are strictly increasing on $\R^+$, locally
 Lipschitz continuous on $\R$ and such that $f(x)>0,g(x)>0$
 for all $x>0$.
For any direction $\nu$ such that $\Lambda_2(\nu)\ne\emptyset$
we have that $\Lambda_0(\nu)\ne\emptyset$. If moreover
$\lambda_0(\nu)<\lambda_2(\nu),$
then there exist a   connected component $C^\nu$ of $\Omega_{\lambda_0(\nu)}
^\nu\setminus Z_{\lambda_0(\nu)}^\nu$
and a connected component $\tilde C^\nu$ of
$\Omega_{\lambda_0(\nu)}^\nu\setminus \tilde
Z_{\lambda_0(\nu)}^\nu$ such that $u\equiv u^\nu_{\lambda_0(\nu)}$
in $C^\nu$ and $v\equiv v^\nu_{\lambda_0(\nu)}$
in $\tilde C^\nu$. For such  components, we have
\begin{equation}
\label{3.6} \nabla u(x)\ne0\quad\forall x\in C^\nu, \quad\nabla
v(x)\ne0\quad\forall
x\in \tilde C^\nu,
\end{equation}
\begin{equation}
\label{3.7}\nabla u(x)=0\quad\forall x\in\partial
C^\nu\setminus(T_{\lambda_0(\nu)}^\nu\cup\partial\Omega),\quad\nabla
v(x)=0\quad\forall x\in\partial
\tilde C^\nu\setminus(T_{\lambda_0(\nu)}^\nu\cup\partial\Omega).
\end{equation}
Moreover, for any $\lambda\in (a(\nu),\lambda_0(\nu))$, we have
\begin{equation}
\label{3.8} u<u_\lambda^\nu\quad in \;\Omega_\lambda^\nu\setminus
Z_\lambda^\nu,\quad v<v_\lambda^\nu\quad in \;\Omega_\lambda^\nu\setminus
\tilde Z_\lambda^\nu,
\end{equation}
and
\begin{equation}
\label{3.9} \frac{\partial u}{\partial\nu}(x)>0\quad\forall
x\in\Omega_{\lambda_0(\nu)}
^\nu\setminus
Z\quad and \quad \frac{\partial v}{\partial\nu}(x)>0\quad\forall
x\in\Omega_{\lambda_0(\nu)}
^\nu\setminus
\tilde Z.
\end{equation}
\end{lemma}
\proof The proof will follow the same steps as in the proof of Theorem 3.1
in \cite{Da}. \\
\mbox{}\\
\underline{\textsc{Step 1:}}  We take a direction $\nu$ such that
$\Lambda_2(\nu)\ne\emptyset$ and we prove that
$\Lambda_0(\nu)\ne\emptyset$. As in the proof
of Theorem 3.1 in \cite{Da},
we can prove
that if $\lambda>a(\nu)$, $\lambda-a(\nu)$ is small, then
$|\Omega_\lambda^\nu|$
is small and since $u\le u_\lambda^\nu$, $v\le v_\lambda^\nu$ on
$\partial\Omega_\lambda^\nu$,
by Theorem \ref{thm2.2}, we get $u\le u_\lambda^\nu$, $v\le v_\lambda^\nu$ on
$\Omega_\lambda^\nu$. So $\Lambda_0(\nu)\ne\emptyset$. Here we have applied
Theorem \ref{thm2.2}
to the pairs $(u,v)$ and $(\bar u,\bar v)=(u_\lambda^\nu,v_\lambda^\nu)$
and with  $A_1=\tilde A_1=\Omega_\lambda^\nu$ and $A_2=\tilde
A_2=A_3=\tilde A_3=\emptyset$.
\\\mbox{}\\
\underline{\textsc{Step 2:}}  By continuity of $u,v$ ,
the inequalities $u\le u_{\lambda_0(\nu)}^\nu$ and $v\le
v_{\lambda_0(\nu)}^\nu$ hold in $\Omega_{\lambda_0(\nu)}^\nu$. Moreover, by
Lemma  \ref{thm1.4}, since $f,g$ are nondecreasing  on $\R^+$,
we have that if $C^\nu$ and $\tilde C^\nu$
are connected components of respectively
$\Omega_{\lambda_0(\nu)}^\nu\setminus Z_{\lambda_0(\nu)}^\nu$
and $\Omega_{\lambda_0(\nu)}^\nu\setminus \tilde
Z_{\lambda_0(\nu)}^\nu$, then either $u<u_{\lambda_0(\nu)}^\nu$ or
$u\equiv u_{\lambda_0(\nu)}^\nu$ in $C^\nu$ and either
$v<v_{\lambda_0(\nu)}^\nu$ or
$v\equiv v_{\lambda_0(\nu)}^\nu$ in $\tilde C^\nu$.
Assume now that $\lambda_0(\nu)<\lambda_2(\nu)$, and by contradiction that
$u<u_{\lambda_0(\nu)}^\nu$ in $\Omega_{\lambda_0(\nu)}^\nu\setminus
Z_{\lambda_0(\nu)}^\nu$. We first show that this implies $v<
v_{\lambda_0(\nu)}^\nu$ on $\Omega_{\lambda_0(\nu)}^\nu\setminus \tilde
Z_{\lambda_0(\nu)}^\nu$. Suppose by contradiction this is not the case and
 that in a
component $\tilde C^\nu$ of $\Omega_{\lambda_0(\nu)}^\nu\setminus \tilde
Z_{\lambda_0(\nu)}^\nu$, we have $v\equiv v_{\lambda_0(\nu)}^\nu$.
 Since $\Omega_{\lambda_0(\nu)}^\nu\setminus \tilde
Z_{\lambda_0(\nu)}^\nu$ is open, $\tilde C^\nu$ is also open.
Moreover,
$int(Z_{\lambda_0(\nu)}^\nu)=\emptyset$. Indeed, if  $Z_{\lambda_0(\nu)}^\nu$
contains an open set $\Omega'$, we have
\begin{equation}
\label{interieur}
\int_{\Omega'}|\nabla u|^{p_1-2}\nabla u.\nabla\phi\,dx=0=
\int_{\Omega'}f(v)\phi\,dx\quad\forall\phi\in
C^\infty_c(\Omega'),
\end{equation}
which is impossible since $v>0$ on $\Omega$ and
$f(x)>0$ for all $x>0$. So
$int(Z_{\lambda_0(\nu)}^\nu)=\emptyset$ and a similar reasoning
would show that $int(\tilde Z_{\lambda_0(\nu)}^\nu)=\emptyset$.
Now if we assume by contradiction that in a
component $\tilde C^\nu$ of $\Omega_{\lambda_0(\nu)}^\nu\setminus \tilde
Z_{\lambda_0(\nu)}^\nu$, we have $v\equiv v_{\lambda_0(\nu)}^\nu$,
then
$u< u_{\lambda_0(\nu)}^\nu$ and $v=v_{\lambda_0(\nu)}^\nu$ in
$\tilde C^\nu\cap(\Omega_{\lambda_0(\nu)}^\nu\setminus
Z_{\lambda_0(\nu)}^\nu)\ne\emptyset$. Denoting $\tilde
C^\nu\cap(\Omega_{\lambda_0(\nu)}^\nu\setminus
Z_{\lambda_0(\nu)}^\nu)$ by $A$, this would imply
$$\int_A|\nabla v|^{p_2-2}\nabla v.\nabla\phi\,dx=\int_A|\nabla
v_{\lambda_0(\nu)}^\nu|^{p_2-
2}\nabla v_{\lambda_0(\nu)}^\nu.\nabla\phi\,dx\quad\forall\phi\in
C^\infty_c(A),$$
and so
$$\int_Ag(u)\phi\,dx=\int_Ag(u_{\lambda_0(\nu)}^\nu)\phi\,dx\quad\forall\phi\in
C^\infty_c(A),$$
which is impossible since $g$ is strictly increasing.
So we also have $v<v_{\lambda_0(\nu)}^\nu$ in
$\Omega_{\lambda_0(\nu)}^\nu\setminus \tilde
Z_{\lambda_0(\nu)}^\nu$.

As in \cite{Da}, we can choose $\alpha,M>0$ independent from
$\lambda\in(a(\nu),\lambda_2(\nu)]$
so that Theorem \ref{thm2.2} applies in $\Omega$ to the pairs $(u,v)$ and
$(\bar u,\bar v)=(u_{\lambda}^\nu,v_{\lambda}^\nu)$ and with
$\Omega'=\Omega_\lambda^\nu$.
Following the same idea as in \cite{Da}, we take two open sets $A,\tilde A$
with $Z_{\lambda_0(\nu)}^\nu\subset A\subset\Omega_{\lambda_0(\nu)}^\nu $
and $\tilde Z_{\lambda_0(\nu)}^\nu\subset \tilde
A\subset\Omega_{\lambda_0(\nu)}^\nu $
such that $\sup_A(|\nabla u|+|\nabla
u_{\lambda_0(\nu)}^\nu|)<\frac{M}{2}$, $\sup_{\tilde A}(|\nabla v|+|\nabla
v_{\lambda_0(\nu)}^\nu|)<\frac{M}{2}$. We also fix a compact set $K\subset
\Omega_{\lambda_0(\nu)}^\nu$ such that $|\Omega_{\lambda_0(\nu)}^\nu
\setminus K|<\frac{\alpha}{2}$. If $K\setminus A\ne\emptyset$,
$\min_{K\setminus
A}(u_{\lambda_0(\nu)}^\nu-u)=m>0$ and if $K\setminus \tilde A\ne\emptyset$,
$\min_{K\setminus
\tilde A}(v_{\lambda_0(\nu)}^\nu-v)=\tilde m>0$. By continuity of $u,v$,
there exists $\epsilon>0$ such that $\lambda_0(\nu)+\epsilon<\lambda_2(\nu)$
and such that for all
$\lambda\in(\lambda_0(\nu),\lambda_0(\nu)+\epsilon)$,
$$|\Omega_{\lambda}^\nu
\setminus K|<\alpha,\quad\sup_{A}(|\nabla u|+|\nabla
u_{\lambda}^\nu|)<M\quad\textrm{and}\quad\sup_{\tilde A}(|\nabla v|+|\nabla
v_{\lambda}^\nu|)<M,$$
$$u_{\lambda}^\nu-u>\frac{m}{2}>0\quad \textrm{in }K\setminus A
\quad\textrm{if}\quad K\setminus  A\ne\emptyset,$$
and
$$v_{\lambda}^\nu-v>\frac{\tilde m}{2}>0\quad \textrm{in }K\setminus
\tilde A\quad\textrm{if}\quad K\setminus \tilde A\ne\emptyset.$$
For such values of $\lambda$, we have that $u\le u_{\lambda}^\nu$
and $v\le v_{\lambda}^\nu$
on respectively $\partial(\Omega_\lambda^\nu\setminus(K\setminus A))$
and $\partial(\Omega_\lambda^\nu\setminus(K\setminus \tilde A))$.
Indeed, if $x_0$ belongs to $\partial(\Omega_\lambda^\nu
\setminus(K\setminus A))$, then either
$x_0\in\partial\Omega_\lambda^\nu$ where trivially $u\le
u_\lambda^\nu$, or $x_0\in\partial(K\setminus A)$, where $u_\lambda^\nu-u$
is positive. The same kind of reasoning works in the case
$x_0\in\partial(\Omega_\lambda^\nu
\setminus(K\setminus \tilde A))$.\\
So we can apply Theorem \ref{thm2.2} to $(u,v)$, $(\bar{u},\bar v)=
(u_\lambda^\nu,v_\lambda^\nu)$, by taking $\Omega'=\Omega_\lambda^\nu$,
$A_1=\tilde A_1=\Omega_\lambda^\nu\setminus
K$, $A_2=K\cap A$, $\tilde A_2= K\cap \tilde A$ and $A_3=K\setminus
A$, $\tilde A_3=K\setminus\tilde A$. We verify easily that $A_1\cup
A_2=\Omega_\lambda^\nu\setminus(K\setminus A)$
and $\tilde A_1\cup \tilde A_2=\Omega_\lambda^\nu\setminus(K\setminus
\tilde A)$
are open subsets contained  in $\Omega_\lambda^\nu$. We
finally conclude that $u\le u_{\lambda}^\nu$
and $v\le v_{\lambda}^\nu$ in $\Omega_\lambda^\nu$ for
$\lambda\in(\lambda_0(\nu),\lambda_0(\nu)+\epsilon)$, which
contradicts the definition of $\lambda_0(\nu)$.\\\mbox{}\\
We prove (\ref{3.6}) and (\ref{3.7}) exactly as in
\cite{Da} (they are simple consequences of the definition of $C^\nu$ and
$\tilde C^\nu$).\\\mbox{}\\
\underline{\textsc{Step 3:}} To prove (\ref{3.8}), it suffices  to
prove
\begin{equation}\label{3.10}
u<u_\lambda^\nu\quad\textrm{in }\Omega^\nu_\lambda\setminus Z\quad\textrm{if
}\lambda\in (a(\nu),\lambda_0(\nu)),
\end{equation}
and
\begin{equation}\label{5.27}
v<v_\lambda^\nu\quad\textrm{in }\Omega^\nu_\lambda\setminus \tilde
Z\quad\textrm{if
}\lambda\in (a(\nu),\lambda_0(\nu)),
\end{equation}
as in \cite{Da}. Indeed,
If (\ref{3.8}) is false, and if (for example)
$u(x_0)=u_\lambda^\nu(x_0)$ for some
point $x_0\in\Omega_\lambda^\nu\setminus Z_\lambda^\nu$, then by
Lemma \ref{thm1.4}, $u=u_\lambda^\nu$ in the component of
$\Omega_\lambda^\nu\setminus
Z_\lambda^\nu$ to which $x_0$
 belongs, and this implies that both $|\nabla u(x_0)|$ and $|\nabla
u_\lambda^\nu(x_0)|$
 are not zero, i.e. $x_0\in\Omega_\lambda^\nu\setminus Z$, so that
 (\ref{3.10}) does not hold. The same reasoning also works for $v$.

Let us prove (\ref{3.10}). The argument is
 the same as in \cite{Da}. We recall  it here for the sake of completeness.
For simplicity of notations, we assume that $\nu=e_1=(1,0,\dots,0)$ and we
omit the
superscript $e_1$ in $\Omega_\lambda^{e_1},\lambda_0(e_1)$,
$a(e_1),u_\lambda^{e_1}$,\dots
We write coordinates in $\R^N$ as $x=(y,z)$ with $y\in\R$,
$z\in\R^{N-1}$.

Arguing by contradiction, we assume that there exists $\mu$ with
$a<\mu<\lambda_0$
and $x_0=(y_0,z_0)\in\Omega_\mu\setminus Z$ such that
$u(x_0)=u_\mu(x_0)$. Since $v\le v_\mu$ on $\Omega_\mu$
and $f$ is nondecreasing,
Lemma \ref{thm1.4} implies $u=u_\mu$ in the component $C$
of $\Omega_\mu\setminus Z_\mu$ to which $x_0$ belongs.
If $\lambda>\mu$, $\lambda-\mu$ is small, $(x_0)_\lambda=x_\mu$
for some point $x=(y,z_0)\in C$ with $y<y_0$. Since for
$\lambda\in(\mu,\lambda_0]$, $u\le u_\lambda$ in $\Omega_\lambda$,
we have $u(y,z_0)=u(x)=u(x_\mu)=u((x_0)_\lambda)\ge
u(x_0)=u(y_0,z_0)$. This implies that $u(y,z_0)=u(y_0,z_0)$ since
by definition of $\lambda_0$, $u$ is nondecreasing in the
$e_1$-direction in $\Omega_{\lambda_0}$. Therefore the
set
$$U:=\{y<y_0\,|\,(y,z_0)\in\Omega\textrm{ and }u(y,z_0)=u(y_0,z_0)\},$$
is not empty. If we set $y_1:=\inf U$, we show that $x_1:=(y_1,z_0)\in\partial
\Omega$. Assume by contradiction that $x_1\in\Omega$ and put
$\lambda_1:=\frac{y_1+y_0}{2}$. By continuity of $u$, $u(x_1)=u(x_0)$
and since $(x_1)_{\lambda_1}=x_0$  and $\nabla u(x_0)\ne0$, we have
$x_1\in\Omega_{\lambda_1}\setminus Z_{\lambda_1}$. By
Lemma \ref{thm1.4}, $u=u_{\lambda_1}$ in the component of
$\Omega_{\lambda_1}\setminus Z_{\lambda_1}$ to which $x_1$
belongs, which implies that $\nabla u(x_1)\ne0$ as above.
Repeating the previous arguments with $\mu,x_0$ substituted by
$\lambda_1,x_1$, we obtain that $u(y,z_0)=u(y_0,z_0)$ for some
$y<y_1$, $y_1-y$ small, and this contradicts the definition of
$y_1$. So $x_1\in\partial\Omega$ and $u(x_1)=0=u(x_0)>0$, a contradiction.
This proves (\ref{3.10}) and hence (\ref{3.8}) for $u$. The
proof of (\ref{5.27}) is similar. \\

The proof of (\ref{3.9}) can be made as in \cite{Da} using the
usual Hopf's lemma.
\cqfd
\mbox{}\\

The following result is the analogue of Proposition 3.1 in
\cite{Da}.
\begin{lemma}\label{prop 3.1}
Suppose that $(u,v)\in C^1_0(\Om)\times C^1_0(\Om) $ is a weak solution of
(\ref{1.27})
where $p_1,p_2\in (1,2)$ and $f,g$ satisfy the hypotheses
of Lemma  \ref{3.1}.
Then for any direction $\nu$, the set
$\Omega_{\lambda_0(\nu)}^\nu $ does not contain any subset $\Gamma$ of $Z$
on which $u$ is constant and whose projection on the hyperplane $T_{\lambda_0
(\nu)}^\nu$ contains a non empty open subset of $T_{\lambda_0
(\nu)}^\nu$ relatively to the induced topology. Similarly, for any
direction $\nu$, the set
$\Omega_{\lambda_0(\nu)}^\nu $ does not contain any subset $\tilde\Gamma$
of $\tilde Z$
on which $v$ is constant and whose projection on the hyperplane $T_{\lambda_0
(\nu)}^\nu$ contains a non empty open subset of $T_{\lambda_0
(\nu)}^\nu$.
\end{lemma}
\proof The proof is identical to that of Proposition 3.1 of \cite{Da}, case
2. We give it here for the sake of completeness. \\

For simplicity of notations, we take $\nu=e_1=(1,0,\dots,0)$ and denote a
point $x\in\R^N$
as $x=(y,z)$ with $y\in\R,z\in\R^{N-1}$. We omit the superscript $\nu=e_1$
in $\Omega_\lambda^\nu,u_\lambda^\nu$, \dots

 Arguing by contradiction, we assume that $\Omega_{\lambda_0}$
 contains a set $\Gamma$ with the properties:
 \begin{itemize}
 \item[(i)] there exists $\gamma>0$ and $z_0\in\R^{N-1}$ such that
 for each $(\lambda_0,z)\in T_{\lambda_0}$ with $|z-z_0|<\gamma$,
 there exists $y<\lambda_0$ with $(y,z)\in\Gamma$,
 \item[(ii)] $\nabla u(x)=0$ for all $x\in\Gamma,$
 \item[(iii)] $u(x)=m>0$ for all $x\in\Gamma$.
 \end{itemize}
 Observe that $\bar\Gamma$ satisfies the same properties as $\Gamma$
 and that by (iii), $\bar\Gamma\cap\partial\Omega=\emptyset$.
 Let $\omega=\omega_\gamma$ be the $(N-1)$-dimensional ball
 centered at $z_0$ with radius $\gamma$. We consider the cylinder
$\R\times\omega$
 and denote by $\Sigma$ the intersection
 $(\R\times\omega)\cap\Omega_{\lambda_0}$. We then define the
``right part of $\Sigma$ with respect to $\Gamma$"
 $$\Sigma_r:=\{(y,z)\in\Sigma\,|\,z\in\omega,\sigma(z)<y<\lambda_0\},$$
 where
 $$\sigma(z):=\sup\{y\in\R\,|\,(y',z)\notin\bar\Gamma\;\textrm{for all
}y'<y\}.$$
 One can see that $\Sigma_r$ is open and well defined and by the
 monotonicity of $u$ in $\Omega_{\lambda_0}$, we have $u(x)\ge m$
 for all $x\in\Sigma_r$. We claim that $u\ne m$ in $\Sigma_r$. Indeed if
 this is not the case i.e. $u=m$
 in $\Sigma_r$, then $\nabla u=0$ in $\Sigma_r$ which is impossible since
 $f(x)>0$ for all $x>0$ and $v>0$ in $\Sigma_r$.
Furthermore,
$$
-\plun(u-m)=-\plun u>0\quad\textrm{on }\Sigma_r.\\
$$
One can see as in \cite{Da} that there exists  $x'$ on
$\partial\Sigma_r\cap\bar\Gamma$
at which
the interior sphere condition is satisfied.
At such a point, we have that $u(x')=m=\min_{\overline{\Sigma_r}}u$.
By the Hopf's lemma, we conclude that $\frac{\partial u}{\partial s}(x')>0$
for
an interior directional derivative, which contradicts (ii) for $\bar \Gamma$.\\

The same arguments can be used to prove the nonexistence of a set
$\tilde\Gamma$.
\cqfd
\mbox{}\\

Before giving the proof of Theorem \ref{mono}, we introduce  some more
notations.
We denote by $\mathcal{F}_\nu$ (resp. $\tilde\mathcal{F}_\nu$)
the collection of the connected components $C^\nu$ of
$\Omega_{\lambda_0(\nu)}^\nu\setminus Z_{\lambda_0(\nu)}^\nu$ (resp.
$\tilde C^\nu$ of
$\Omega_{\lambda_0(\nu)}^\nu\setminus\tilde Z_{\lambda_0(\nu)}^\nu$) such
that $u=u_{\lambda_0(\nu)}^\nu$ in $C^\nu$, $\nabla u\ne0$ in $C^\nu$
and $\nabla u=0$ on $\partial C^\nu\setminus(T_{\lambda_0(\nu)}^\nu\cup
\partial\Omega)$ (resp. $v=v_{\lambda_0(\nu)}^\nu$ in $\tilde C^\nu$,
$\nabla v\ne0$ in $\tilde C^\nu$
and $\nabla v=0$ on $\partial \tilde C^\nu\setminus( T_{\lambda_0(\nu)}^\nu\cup
\partial\Omega)$ ). Then we define
\begin{equation}
\label{CCprime}C'^\nu:=C^\nu,
\end{equation}
 if $\nabla
u(x)=0$ for all $x\in \partial C^\nu$, or, if there are some points
\label{111} $x\in\partial C^\nu\cap T_{\lambda_0(\nu)}^\nu$ such that $\nabla
u(x)\ne0$, we define
 \begin{equation}\label{Cprime}C'^\nu:=C^\nu\cup C^\nu_1\cup
C^\nu_2,\end{equation} where
 $C_1^\nu=R_{\lambda_0(\nu)}(C^\nu)$, $C_2^\nu=\{x\in\partial
 C^\nu\cap T_{\lambda_0(\nu)}\,|\,\nabla u(x)\ne0\}$. As in \cite{Da}, we
 can check that $C'^\nu$ is open and connected, with
\begin{equation}\label{propr}\nabla u\ne0\quad\textrm{in }
C'^\nu, \qquad\nabla u=0\quad\textrm{on }\partial C'^\nu.
\end{equation}
We define in the same way
$\tilde C'^\nu$ with $u$ replaced by $v$. Finally we define
$\mathcal{F}'_\nu:=\{C'^\nu\,|\,C^\nu\in\mathcal{F}_\nu\}$ and
$\tilde \mathcal{F}'_\nu:=\{\tilde C'^\nu\,|\,\tilde
C^\nu\in\tilde\mathcal{F}_\nu\}$.\\
\begin{remark}
\rm\label{rem4.1}We have the analogue of remark 4.1 from \cite{Da} for $u$
and $v$: if
$\nu_1,\nu_2$ are two directions, then either $C'^{\nu_1}=C'^{\nu_2}$ or
$C'^{\nu_1}\cap C'^{\nu_2}=\emptyset$. Indeed,
if $C^{\nu_1}\in\mathcal{F}_{\nu_1},
C^{\nu_2}\in\mathcal{F}_{\nu_2}$, if $C'^{\nu_1}\cap
C'^{\nu_2}\ne\emptyset$ and $C'^{\nu_1}\ne
C'^{\nu_2}$ , then $\partial C'^{\nu_1}\cap
 C'^{\nu_2}\ne\emptyset$ or $C'^{\nu_1}\cap
 \partial C'^{\nu_2}\ne \emptyset$ by Corollary 4.1 of \cite{Da}, which is
 impossible since $\nabla u\ne0$ in $C'^{\nu_i}$, $\nabla u=0$ on
 $\partial C'^{\nu_i}$ for $i=1,2$. So necessarily, either $C'^{\nu_1}=
C'^{\nu_2}$ or $C'^{\nu_1}\cap
C'^{\nu_2}=\emptyset$. In a similar way, if $\tilde C^{\nu_1}\in\tilde
\mathcal{F}_{\nu_1},
\tilde C^{\nu_2}\in\tilde \mathcal{F}_{\nu_2}$, then either $\tilde C'^{\nu_1}=
\tilde C'^{\nu_2}$ or $\tilde C'^{\nu_1}\cap
\tilde C'^{\nu_2}=\emptyset$.
\end{remark}

 \mbox{}\\
\textbf{Proof of Theorem \ref{mono}:}
The proof will follow the steps of that of Theorem 1.1 from \cite{Da}. \\
If $\nu $ is a direction and $\delta>0$, we denote  by $I_\delta(\nu)$ the set
$$
I_\delta(\nu):=\{\mu\in\R^N\,|\,|\mu|=1,|\mu-\nu|<\delta\}.
$$
 As in the proof of Lemma \ref{3.1}, we can fix $\alpha,M>0$ such that
Theorem \ref{thm2.2} applies to any direction $\nu$ and any set
$\Omega'\subset\Omega_\lambda^\nu$
for all $\lambda\in(a(\nu),\lambda_1(\nu)]$.\\
Suppose by contradiction that $\nu_0$ is a direction such that
$\lambda_0(\nu_0)<\lambda_1(\nu_0)$. Since $\lambda_1(\nu_0)\le
\lambda_2(\nu_0)$, then $\lambda_0(\nu_0)<\lambda_2(\nu_0)$ and by
Lemma \ref{3.1}, $\mathcal{F}_{\nu_0}$ and  $\tilde\mathcal{F}_{\nu_0}$
are non empty. Since $\R^N$ is separable and since each component of
$\mathcal{F}_{\nu_0}$ and  $\tilde \mathcal{F}_{\nu_0}$ is open,
$\mathcal{F}_{\nu_0}$
and  $\tilde \mathcal{F}_{\nu_0}$ contain at most countable many sets,
so  $\mathcal{F}_{\nu_0}=\{C_i^{\nu_0}\,|\,i\in I\subseteq\N\}$
and  $\tilde\mathcal{F}_{\nu_0}=\{\tilde C_i^{\nu_0}\,|\,i\in \tilde I
\subseteq\N\}$. In case  $I$  is infinite, since the components
are disjoint, we have that $\sum_{i=1}^\infty|C_i^{\nu_0}|\le|\Omega|$, so
we can choose $n_0\ge1$ for which
$$\sum_{i=n_0+1}^\infty|C_i^{\nu_0}|<\frac{\alpha}{6},$$
 and the same remark holds for $\tilde I$ with a number $\tilde n_0$. If $I$
and $\tilde I$ are finite, let $n_0$ and $\tilde n_0$ be their cardinality.
We then choose two compacts
$K_0\subset(\Omega_{\lambda_0(\nu_0)}^{\nu_0}\setminus \cup_{i\in
I}C_i^{\nu_0})$ and $\tilde
K_0\subset(\Omega_{\lambda_0(\nu_0)}^{\nu_0}\setminus
 \cup_{i\in
\tilde I}\tilde C_i^{\nu_0})$ such that
$$
|(\Omega_{\lambda_0(\nu_0)}^{\nu_0}\setminus \cup_{i\in
I}C_i^{\nu_0})\setminus K_0|<\frac{\alpha}{6}\quad\textrm{and}\quad
|(\Omega_{\lambda_0(\nu_0)}^{\nu_0}\setminus \cup_{i\in
\tilde I}\tilde C_i^{\nu_0})\setminus \tilde K_0|<\frac{\alpha}{6}.$$
Then we take some compact sets $K_i\subset C_i^{\nu_0}, i=1,\dots, n_0$
and  $\tilde K_i\subset \tilde C_i^{\nu_0}, i=1,\dots, \tilde n_0$, such
that
$$
|C_i^{\nu_0}\setminus  K_i|<\frac{\alpha}{6n_0}\quad i=1,\dots,
n_0\quad\textrm{and}\quad
|\tilde C_i^{\nu_0}\setminus \tilde K_i|<\frac{\alpha}{6\tilde n_0}\quad
i=1,\dots,
\tilde n_0.$$
So we have decomposed $\Omega_{\lambda_0(\nu_0)}^{\nu_0}$ in the sets
$K_0,K_1,\dots,K_{n_0}$ and in a remaining part with measure less than
$\frac{\alpha}{2}$, and the same remark holds with the sets $\tilde K_i$.\\
We define then $A:=\{x\in\Omega_{\lambda_0(\nu_0)}^{\nu_0}\,|\,|\nabla
u(x)|+|\nabla u_{\lambda_0(\nu_0)}^{\nu_0}(x)|<\frac{M}{2}\}$ and
$\tilde A:=\{x\in\Omega_{\lambda_0(\nu_0)}^{\nu_0}\,|\,|\nabla
v(x)|+|\nabla v_{\lambda_0(\nu_0)}^{\nu_0}(x)|<\frac{M}{2}\}$. Clearly, the
sets
$K_0\setminus A$ and $\tilde K_0\setminus \tilde A$ are compact.
By
 Lemma
\ref{thm1.4} and since $K_0\subset
(\Omega_{\lambda_0(\nu_0)}^{\nu_0}\setminus \cup_{i\in
I}C_i^{\nu_0})$ and $\tilde
K_0\subset(\Omega_{\lambda_0(\nu_0)}^{\nu_0}\setminus
 \cup_{i\in
\tilde I}\tilde C_i^{\nu_0})$, we have that  $u<
u_{\lambda_0(\nu_0)}^{\nu_0}$ in $K_0\setminus A$ if $K_0\setminus
A\ne\emptyset$ and
$v<v_{\lambda_0(\nu_0)}^{\nu_0}$ in $\tilde K_0\setminus \tilde A$ if $\tilde
K_0\setminus \tilde A\ne\emptyset$. So  there exists $m>0$ such that
$$
u_{\lambda_0(\nu_0)}^{\nu_0}-u\ge m>0\quad\textrm{in }K_0\setminus A
\quad\textrm{if}\quad K_0\setminus A\ne\emptyset,$$and$$
v_{\lambda_0(\nu_0)}^{\nu_0}-v\ge m>0\quad\textrm{in }\tilde K_0\setminus
\tilde A\quad\textrm{if}\quad K_0\setminus\tilde A\ne\emptyset.
$$
Since $\Omega$ is of class $C^1$, $a(\nu)$ is continuous and
$\lambda_1(\nu)$ is lower semicontinuous
(see \cite{Az2}). By continuity, there exists
$\epsilon_0,\delta_0>0$ such that if
$|\lambda-\lambda_0(\nu_0)|\le\epsilon_0$ and $|\nu-\nu_0|\le\delta_0$,
then
\begin{equation}\label{4.1}
\lambda_1(\nu)>\lambda_0(\nu_0)+\epsilon_0,
\end{equation}
\begin{equation}\label{4.2}
|\nabla u|+|\nabla u_\lambda^\nu|<M\textrm{ on }A\quad\textrm{and}\quad
 |\nabla v|+|\nabla v_\lambda^\nu|<M\textrm{ on } \tilde A,
\end{equation}
$$\begin{array}{rcl}
K_i\subset\Omega_\lambda^\nu\quad \textrm{for
}i=0,\dots,n_0&\textrm{and}&
\tilde K_i\subset\Omega_\lambda^\nu\quad \textrm{for }i=0,\dots,\tilde
n_0,\vspace*{0.2cm}\\
R_\lambda^\nu(K_i)\subset R_\lambda^\nu(C_i^{\nu_0})\quad \textrm{for
}i=1,\dots,n_0&\textrm{and}&
 R_\lambda^\nu(\tilde K_i)\subset R_\lambda^\nu(\tilde C_i^{\nu_0})\quad
\textrm{for }i=1,\dots,\tilde n_0,
\end{array}$$
$$
|\Omega_\lambda^\nu\setminus\cup_{i=0}^{n_0}K_i|<\alpha\quad\textrm{and}\quad
|\Omega_\lambda^\nu\setminus\cup_{i=0}^{\tilde n_0}\tilde K_i|<\alpha
$$
and
\begin{equation}\label{4.5}
u_\lambda^\nu-u\ge\frac{m}{2}\quad\textrm{in }K_0\setminus A\quad
\textrm{and}\quad v_\lambda^\nu-v\ge\frac{m}{2}\quad\textrm{in }\tilde
K_0\setminus
\tilde A.
\end{equation}

We now proceed in several steps in order to show that there exists
$i_1\in\{1,\dots, n_0\}$ or $j_1\in\{1,\dots,\tilde n_0\}$ and a
direction $\nu_1\in I_{\delta_0}$ such that
$C'^{\nu_0}_{i_1}\in\mathcal{F}'_\nu$
for any direction $\nu$ in a suitable neighborhood $I_\delta(\nu_1)$
of $\nu_1$  and $\partial C^{\nu_0}_{i_1}$ contains a set $\Gamma$
as in Lemma \ref{prop 3.1} (with respect to the direction $\nu_1$)
or such that $\tilde C'^{\nu_0}_{j_1}\in\tilde\mathcal{F}'_\nu$
for any direction $\nu$ in a suitable neighborhood $I_\delta(\nu_1)$
of $\nu_1$ and $\partial \tilde C^{\nu_0}_{j_1}$ contains a set $\tilde\Gamma$
as in Lemma \ref{prop 3.1}. This together with Lemma \ref{prop 3.1} would
lead to a contradiction and end the proof.\\\mbox{}\\
\underline{\textsc{Step 1:}} We show that $\lambda_0(\nu)$ is continuous
at $\nu_0$ with respect to
$\nu$ i.e. for
all $\epsilon\in(0,\epsilon_0)$, there exists $\delta\in(0,\delta_0)$ such that
if $\nu\in I_{\delta}(\nu_0)$, then
$$
(i) \quad
\lambda_0(\nu_0)-\epsilon<\lambda_0(\nu)<\lambda_0(\nu_0)+\epsilon.
$$
Moreover, for all $\nu\in I_{\delta(\epsilon_0)}(\nu_0)$, we have
$$
(ii)\quad \exists i\in\{1,\dots,n_0\}\textrm{ such that }C
_i^{'\nu_0}\in\mathcal{F}'_\nu\quad \textrm{or}\quad \exists
i\in\{1,\dots,\tilde n_0\}\textrm{ such that }\tilde C
_i^{'\nu_0}\in\tilde \mathcal{F}'_\nu.
$$
\textsc{Proof of Step 1:} Let $0<\epsilon\le\epsilon_0$ be fixed. By
definition of
$\lambda_0(\nu_0)$, there exists
$\lambda\in(\lambda_0(\nu_0),\lambda_0(\nu_0)+\epsilon)$ and
$x\in\Omega_\lambda^{\nu_0}$ such that $u(x)>u_\lambda^{\nu_0}(x)$
or $\tilde x\in\Omega_\lambda^{\nu_0}$ such that
$v(\tilde x)>v_\lambda^{\nu_0}(\tilde x)$.
Suppose we are in the first case. Then by continuity, there
exists $\delta_1\in(0,\delta_0]$ such that for all $\nu\in
I_{\delta_1}(\nu_0)$, $x$
 belongs to $\Omega_\lambda^{\nu}$ and
$u(x)>u_\lambda^\nu(x)$.
This implies that for all $\nu\in I_{\delta_1}(\nu_0)$, we have
$\lambda_0(\nu)<\lambda<\lambda_0(\nu_0)+\epsilon$. We conclude in the same way
if we are in the second case.\\

We then show that
there exists $\delta_2\in(0,\delta_0]$ such that
$\lambda_0(\nu)>\lambda_0(\nu_0)-\epsilon$ for all $\nu\in
I_{\delta_2}(\nu_0)$. Suppose that this is false. Thus we can find a
sequence $\nu_n\to\nu_0$ with $\lambda_0(\nu_n)\le\lambda_0(\nu_0)-\epsilon$
for all $n\in\N_0$. Passing to a subsequence still denoted by $\nu_n$, we have
$\lambda_0(\nu_n)\to\bar\lambda\le\lambda_0(\nu_0)-\epsilon$.
Since $\lambda_0(\nu_n)>a(\nu_n)$ and $a(\nu_n)\to a(\nu_0)$ by the
smoothness of
$\partial\Omega$ (see \cite{Az2}), we
also have $a(\nu_0)\le\bar\lambda$. In fact, this inequality is
strict since by
Theorem \ref{thm2.2} and by definition of $\lambda_0(\nu_n)$,
$|\Omega_{\lambda_0(\nu_n)}^{\nu_n}|\ge\alpha$  thus we  also have
$|\Omega_{\bar\lambda}^{\nu_0}|>0$ which implies
$\bar{\lambda}>a(\nu_0)$. Since $\bar\lambda<\lambda_0(\nu_0)$, by
(\ref{3.8}) of Lemma \ref{3.1}, we have
\begin{eqnarray}
\label{4.6}
u<u_{\bar\lambda}^{\nu_0}&\textrm{ in
}&\Omega_{\bar\lambda}^{\nu_0}\setminus Z_{\bar\lambda}^{\nu_0},\\
\label{4.6'}v<v_{\bar\lambda}^{\nu_0}&\textrm{ in
}&\Omega_{\bar\lambda}^{\nu_0}\setminus \tilde Z_{\bar\lambda}^{\nu_0}.
\end{eqnarray}
Now we make the same reasoning as in step 2 of the proof of Lemma \ref{3.1}.
We can construct  open sets $A,\tilde A\subset\Omega_{\bar\lambda}^{\nu_0}$
and a compact set $K\subset\Omega_{\bar\lambda}^{\nu_0}$ such that
$$
Z_{\bar\lambda}^{\nu_0}\subset A,\quad \sup_A(|\nabla u|+|\nabla
u_{\bar\lambda}^{\nu_0}|)<\frac{M}{2},\quad|\Omega_{\bar\lambda}^{\nu_0}
\setminus K|<\frac{\alpha}{2},
$$
$$
\tilde Z_{\bar\lambda}^{\nu_0}\subset \tilde A,\quad \sup_{\tilde A}
(|\nabla v|+|\nabla
v_{\bar\lambda}^{\nu_0}|)<\frac{M}{2},
$$
and
\begin{eqnarray*}
u_{\bar\lambda}^{\nu_0}-u\ge \gamma>0&\textrm{in}&K\setminus A\quad
\textrm{if}\quad K\setminus A\ne\emptyset,\\
v_{\bar\lambda}^{\nu_0}-v\ge \gamma>0&\textrm{in}&K\setminus \tilde A\quad
\textrm{if}\quad K\setminus \tilde
A\ne\emptyset,
\end{eqnarray*}
for some $\gamma>0$. By continuity of $u,v$ and their gradients, there
exist $r,\delta>0$ such that
$$\sup_A(|\nabla u|+|\nabla
u_{\lambda}^{\nu}|)<M,\quad|\Omega_{\lambda}^{\nu}\setminus
K|<\alpha,\quad u_{\lambda}^{\nu}-u\ge\frac{\gamma}{2}>0\quad\textrm{in
}K\setminus A,
$$
$$\sup_{\tilde A}(|\nabla v|+|\nabla
v_{\lambda}^{\nu}|)<M,\quad
v_{\lambda}^{\nu}-v\ge\frac{\gamma}{2}>0\quad\textrm{in
}K\setminus \tilde A
$$
for all $\nu\in I_r(\nu_0)$ and
$\lambda\in(\bar\lambda-\delta,\bar\lambda+\delta)$. We then apply
Theorem \ref{thm2.2} to the pairs $(u,v)$ and $(\bar u,\bar v)=(u_\lambda^\nu,
v_\lambda^\nu)$, with $\Omega'=\Omega_\lambda^\nu, A_1=\tilde
A_1=\Omega_\lambda^\nu\setminus K$, $A_2=K\cap A,\tilde
A_2=K\cap\tilde A$, $A_3=K\setminus A, \tilde A_3=K\setminus \tilde
A$ and we obtain
$$u\le u_\lambda^\nu\quad \textrm{and}\quad v\le v_\lambda^\nu
\quad\textrm{on}\;
\Omega_\lambda^\nu \quad\forall\nu\in I_r(\nu_0),\,\forall\lambda\in(\bar
\lambda-\delta,\bar\lambda+\delta).
$$ Taking
$\nu=\nu_n,\lambda=\lambda_0(\nu_n)+\eta$ for
$n$ large and $\eta>0$ small, this contradicts the definition of
$\lambda_0(\nu_n)$ and proves (i). Notice that everything above has some sense
since $\lambda<\lambda_1(\nu)$ for $\lambda$ close to $\bar\lambda$ and $\nu$
close to $\nu_0$ by the lower semicontinuity of $\lambda_1(\nu)$ and since
$\bar\lambda<\lambda_1(\nu_0)-\epsilon$. Remark also that the proof of (i)
does not
use the fact that $\lambda_0(\nu_0)<\lambda_2(\nu_0)$, even when we use Lemma
\ref{3.1}, so that the continuity of $\lambda_0(\nu)$ is in fact insured at
all $\nu$
not necessarily equal to $\nu_0$.\\

Observe that since $\epsilon\le\epsilon_0$ and
$\delta\le\delta_0$, by (i) and (\ref{4.1}), we have
$\lambda_1(\nu)>\lambda_0(\nu)$
for all $\nu\in I_{\delta(\epsilon_0)}(\nu_0)$ and
(\ref{4.2})--(\ref{4.5}) are true for $\nu\in
I_{\delta(\epsilon_0)}(\nu_0)$, $\lambda=\lambda_0(\nu)$.\\

Let us now prove (ii). We fix a direction $\nu\in
I_{\delta(\epsilon_0)}(\nu_0)$. Suppose that there exist $i\in\{1,\dots,n_0\}$
and a point $x_i\in K_i$ such that
$u(x_i)=u_{\lambda_0(\nu)}^\nu(x_i)$. Since $\nabla u(x_i)\ne 0$,
by Lemma \ref{thm1.4}, we have $u=u_{\lambda_0(\nu)}^\nu$ in
the component $C^\nu$ of $\mathcal{F}_\nu$ to which $x_i$ belongs.
Since $K_i\subset C_i^{\nu_0}\subset
C'^{\nu_0}_i$ it follows that $x_i\in C'^{\nu}\cap C_i^{'\nu_0}$, hence by
the analogue of remark 4.1 of \cite{Da}, it holds that
$C'^{\nu}=C'^{\nu_0}$. As a consequence, we conclude that
(ii) holds with an index $i\in\{1,\dots,n_0\}$. The same argument works for
$v$: if
there exists $\tilde x_i\in \tilde K_i$ for some $i\in\{1,\dots,\tilde n_0\}$
such that $v(\tilde x_i)=v_{\lambda_0(\nu)}^\nu(\tilde x_i)$, then
$\tilde C'^{\nu_0}_i\in\tilde \mathcal{F}_\nu'$.\\

Next we analyze the case when $u<u_{\lambda_0(\nu)}^\nu$
on $\cup_{i=1}^{n_0}K_i$ and $v<v_{\lambda_0(\nu)}^\nu$
on $\cup_{i=1}^{n_0}\tilde K_i$. If we take
$\lambda>\lambda_0(\nu)$, $|\lambda-\lambda_0(\nu)|$ small, then
by (\ref{4.2})--(\ref{4.5}),
we also have $u<u_{\lambda}^\nu$ on $\cup_{i=1}^{n_0}K_i\cup (K_0\setminus A)$,
$v<v_{\lambda}^\nu$ on $\cup_{i=1}^{\tilde n_0}\tilde K_i\cup(\tilde K_0
\setminus\tilde
A)$,
$\sup_{K_0\cap A}(|\nabla u|+|\nabla u_{\lambda}^\nu|)<M$ ,
$\sup_{\tilde K_0\cap \tilde A}|\nabla v|+|\nabla v_{\lambda}^\nu|<M$,
and finally $|\Omega_\lambda^\nu\setminus\cup_{i=0}^{n_0}K_i|<\alpha$
and $|\Omega_\lambda^\nu\setminus\cup_{i=0}^{n_0}\tilde
K_i|<\alpha$. Applying Theorem \ref{thm2.2} to the pairs  $(u,v)$, $(\bar
u,\bar v)=
(u_\lambda^\nu,v_\lambda^\nu)$, with $\Omega'=\Omega_\lambda^\nu$ and
$$
A_1=\Omega_\lambda^\nu\setminus\cup_{i=0}^{n_0}K_i,\quad \tilde
A_1=\Omega_\lambda^\nu\setminus\cup_{i=0}^{\tilde n_0}\tilde
K_i,
$$
$$
A_2=K_0\cap A,\quad\tilde A_2=\tilde K_0\cap\tilde A,
$$
$$
A_3=\cup_{i=1}^{n_0}K_i\cup (K_0\setminus A),\quad\tilde A_3=\cup_{i=1}^
{\tilde n_0}\tilde K_i\cup(\tilde K_0\setminus\tilde
A),
$$
we get $u\le u_\lambda^\nu$ and $v\le v_\lambda^\nu$ on
$\Omega_\lambda^\nu$, which contradicts the definition of
$\lambda_0(\nu)$. So we conclude that the case $u<u_{\lambda_0(\nu)}^\nu$
on $\cup_{i=1}^{n_0}K_i$ and $v<v_{\lambda_0(\nu)}^\nu$
on $\cup_{i=1}^{n_0}\tilde K_i$ is impossible. This proves (ii).\\

\underline{\textsc{Step 2:}} We prove that there exist $\nu_1\in
I_{\delta(\epsilon_0)}(\nu_0)$ and a neighborhood $I_{\delta_1}(\nu_1)$
such that either there exists an index $i_1\in\{1,\dots,n_0\}$ such that
$$C'^{\nu_0}_{i_1}\in\mathcal{F}'_\nu \quad\forall\nu\in
I_{\delta_1}(\nu_1),
$$
or an index $j_1\in\{1,\dots,\tilde n_0\}$ such that
$$
\tilde C'^{\nu_0}_{j_1}\in\tilde
\mathcal{F}'_\nu\quad\forall\nu\in
I_{\delta_1}(\nu_1).$$

Observe that in the proof
of (ii) of step 1, we have seen that if $\nu\in I_{\delta(\epsilon_0)}(\nu_0)$
and if there exists $x_i\in K_i$ for some $i\in\{1,\dots,n_0\}$
such that $u(x_i)=u_{\lambda_0(\nu_0)}^{\nu_0}(x_i)$, then $u\equiv
u^{\nu_0}_{\lambda_0(\nu_0)}\equiv u^{\nu}_{\lambda_0(\nu)}$ in $K_i$
and $C'^{\nu_0}_i\in\mathcal{F}'_\nu$. So for any $\nu\in
I_{\delta(\epsilon_0)}(\nu_0)$
and each $i\in\{1,\dots,n_0\}$, we have the alternative: either $u\equiv
u^{\nu}_{\lambda_0(\nu)}$
in $K_i$ and $C'^{\nu_0}_i\in\mathcal{F}'_\nu$, or $u<u_{\lambda_0(\nu)}^\nu$
in $K_i$ (by Lemma \ref{thm1.4}, using the monotonicity of $f$ and the
definition of $\lambda_0(\nu)$). In the
latter case, since $K_i$ is compact, $u_{\lambda_0(\nu)}^\nu-u\ge m>0$
in $K_i$ for some $m$, and since the function $\nu\mapsto\lambda_0(\nu)$
is continuous by step 1, we get that the inequality
$u<u_{\lambda_0(\mu)}^\mu$ holds in $K_i$ for all $\mu$ in a
suitable neighborhood $I(\nu)$ of $\nu$ and so $C'^{\nu_0}_i\notin\mathcal{F}'
_\mu$
for all $\mu\in I(\nu)$.
The same remark holds for $v$ and any component $\tilde
C'^{\nu_0}_{j}$ for $j\in\{1,\dots,\tilde n_{0}\}$.

 We start  by taking two sets in
$\mathcal{F}'_{\nu_0}$ and $\tilde \mathcal{F}'_{\nu_0}$, say $C'^{\nu_0}_1$
and $\tilde C'^{\nu_0}_1$ and we argue as
follows. If $C'^{\nu_0}_1\in\mathcal{F}'_\nu$ for any $\nu\in
I_{\delta(\epsilon_0)}(\nu_0)$ or  $\tilde C'^{\nu_0}_1\in\tilde
\mathcal{F}'_\nu$ for any $\nu\in
I_{\delta(\epsilon_0)}(\nu_0)$, then the assertion is proved.
Otherwise, for what we explained above, for some $\mu_1\in
I_{\delta(\epsilon_0)}(\nu_0)$, we have
$C'^{\nu_0}_1\notin\mathcal{F}'_{\mu_1}$
and $u<u_{\lambda_0(\mu_1)}^{\mu_1}$ in $K_1$, and so $C'^{\nu_0}_1\notin
\mathcal{F}'_\mu$ for any $\mu$ in a suitable neighborhood
$I_{\delta_1}(\mu_1)$
of $\mu_1$
which can be chosen such that $I_{\delta_1}(\mu_1)\subset
I_{\delta(\epsilon_0)}(\nu_0)$. If $n_{0}=\tilde n_{0}=1$, then by
(ii) of step 1, $\tilde C'^{\nu_0}_1\in\tilde
\mathcal{F}'_\mu$ for any $\mu$ in $I_{\delta_1}(\mu_1)$
 and the assertion is proved. So we can suppose
$n_{0}+\tilde n_{0}\ge 3$.
Now, by (ii) of step 1, there
exists $i\in\{2,\dots,n_{0}\}$
such that $C'^{\nu_0}_i\in\mathcal{F}'_{\mu_1}$ or
  $j\in\{1,\dots,\tilde n_{0}\}$ such that
 $\tilde C'^{\nu_0}_j\in\tilde\mathcal{F}'_{\mu_1}$. Let us denote
 by $A_2$
 this set $C'^{\nu_0}_i$ or $\tilde C'^{\nu_0}_j$ and suppose
 $A_2\in \mathcal{F}'_{\mu_1}$ (the
 argument is similar in the other case $A_2\in
\tilde\mathcal{F}'_{\mu_1}$). If $A_{2}\in\mathcal{F}'_{\mu} $
 for all $\mu\in I_{\delta_1}(\mu_1)$, then the assertion is proved.
 Otherwise, there exists $\mu_{2}\in I_{\delta_1}(\mu_1)$ such that
 $A_{2}\notin\mathcal{F}'_{\mu_{2}}$ and as above, this
 implies that $A_{2}\notin\mathcal{F}'_{\mu}$ for all $\mu $ in a
 suitable neighborhood  $I_{\delta_2}(\mu_2)$ which can be chosen
 such that $I_{\delta_2}(\mu_2)\subset I_{\delta_1}(\mu_1)$. If
 $n_{0}+\tilde n_{0}=3$, by (ii) of step 1, this implies that there exists
a set $A_{3}$
 in $\mathcal{F}'_{\nu_0}$ or $\tilde \mathcal{F}'_{\nu_0}$
 ($A_3=\tilde C'^{\nu_0}_1$ if $A_2\in\mathcal{F}'_{\nu_0}$)
 with $A_{3}\in \mathcal{F}'_{\mu}$ for all $\mu $
 in  $I_{\delta_2}(\mu_2)$ (or $A_{3}\in \tilde \mathcal{F}'_{\mu}$ for all
$\mu $
 in  $I_{\delta_2}(\mu_2)$) and the assertion is proved. If $n_{0}+\tilde
 n_{0}>3$, we proceed as before taking a set $A_{3}$ in  $\mathcal{F}'_{\nu_0}$
 (or $\tilde \mathcal{F}'_{\nu_0}$) such that $A_{3}\in
 \mathcal{F}'_{\mu_{2}}$ (or in  $ \tilde\mathcal{F}'_{\mu_{2}}$).
 Arguing as we did before, after $k<n_0+\tilde n_{0}$ steps we get a
 set $A_{k}$ in $\mathcal{F}'_{\nu_0}$
 (or $\tilde \mathcal{F}'_{\nu_0}$) with $A_{k}\in\mathcal{F}'_{\mu}$ for
all $\mu $
 in  $I_{\delta_{k-1}}(\mu_{k-1})$ (or $A_{k}\in\tilde\mathcal{F}'_{\mu}$
 for all $\mu $
 in  $I_{\delta_{k-1}}(\mu_{k-1})$) proving the assertion,
 or after $n_0+\tilde n_{0}$ steps, we get a direction
 $\mu_{n_{0}+\tilde n_{0}}\in I_{\delta(\epsilon_0)}(\nu_0)$ such
 that $C'^{\nu_0}_i\notin \mathcal{F}'_{\mu_{n_0+\tilde n_0}}$ for all
$i\in\{1,\dots, n_{0}\}$
 and
 $\tilde C'^{\nu_0}_j\notin \tilde \mathcal{F}'_{\mu_{n_0+\tilde n_0}}$
 for all
 $j\in\{1,\dots,\tilde n_{0}\}$, a contradiction with (ii) of step
 1.\\

\mbox{}\\
\underline{\textsc{Step 3:}} Here we prove that if step 2
is satisfied with a set $C_1'^{\nu_0}\in\mathcal{F}'_{\nu_0}$,
then the set
$\Omega_{\lambda_0(\nu)}^\nu $  contains a subset $\Gamma$ of $Z$
on which $u$ is constant and whose projection on the hyperplane $T_{\lambda_0
(\nu)}^\nu$ contains an open subset of $T_{\lambda_0
(\nu)}^\nu$ relatively to the induced topology. In the same way,  if  step
2
is satisfied with a set $\tilde C'^{\nu_0}_i\in\tilde\mathcal{F}'_{\nu_0}$,
then $\Omega_{\lambda_0(\nu)}^\nu $
contains a subset $\tilde \Gamma$ of $\tilde Z$
on which $v$ is constant and whose projection on the hyperplane $T_{\lambda_0
(\nu)}^\nu$ contains an open subset of $T_{\lambda_0
(\nu)}^\nu$. In both cases, we arrive to a contradiction with
Lemma \ref{prop 3.1} and Theorem \ref{mono} is proved. We shall proceed
exactly as
 in   step 3 of the proof of Theorem 1.1 from
\cite{Da}. We recall the argument for the sake of completeness. \\

Suppose  step 2 is satisfied with some set
$C'^{\nu_0}_i\in\mathcal{F}'_{\nu_0}$
(the same reasoning can be made in the other case,
replacing $u$ by $v$ in what follows).
Let us assume that  $\nu_1,i_1,\delta_1$ are as in step 2 and denote
by $C$ the set $C_{i_1}^{\nu_0}$ from step 2.  We assume moreover for
simplicity that $\nu_1=e_1=(1,0,\dots,0)$ and we omit $\nu_1$ in
$\lambda_0(\nu_1), \Omega_{\lambda_0(\nu_1)}^{\nu_1},
T_{\lambda_0(\nu_1)}^{\nu_1}$,\dots

By step 2 and definition of $C'$ (see (\ref{CCprime},\ref{Cprime})), we
have that for
each $\nu\in I_{\delta_1}(e_1)$,
\begin{equation}\label{4.8}
u\equiv u_{\lambda_0(\nu)}^\nu\quad\textrm{in }C'\quad\textrm{and} \quad
C'\in\mathcal{F}'_\nu.
\end{equation}
Since $C'$ is open (see the definition of $C'$),
$R_{\lambda_0(\nu)}^\nu(C')\cap
R_{\lambda_0}(C')\ne\emptyset$ if $\nu$ is sufficiently close to $e_1$.
Moreover, by
(\ref{propr}), $\nabla u=0$ on $\partial C'$ which implies by (\ref{4.8})
that
\begin{equation}\label{part}
\nabla u=0\textrm{ on }\partial R_{\lambda_0(\nu)}^\nu(C')\quad
\textrm{for all }\nu\in
I_{\delta_1}(e_1).
\end{equation}
Arguing as in Remark \ref{rem4.1}, if $R_{\lambda_0(\nu)}^\nu(C')\cap
R_{\lambda_0}(C')\ne\emptyset$ and $R_{\lambda_0(\nu)}^\nu(C')\ne
R_{\lambda_0}(C')$, then by Corollary 4.1 from \cite{Da},
$\partial(R_{\lambda_0(\nu)}^\nu(C'))\cap
R_{\lambda_0}(C')\ne\emptyset$ or $R_{\lambda_0(\nu)}^\nu(C')\cap\partial(
R_{\lambda_0}(C'))\ne\emptyset$, which is impossible by
(\ref{propr}), (\ref{4.8}) and (\ref{part}). Thus, we get that
$$R_{\lambda_0(\nu)}^\nu(C')\equiv R_{\lambda_0}(C')$$
for $\nu$
sufficiently close to $e_1$, say for $\nu\in I_{\delta_2}(e_1)$
for some $0<\delta_2\le\delta_1$. Then we take a point $\bar x=(\bar y,\bar z)$
in $\partial C\cap \Omega_{\lambda_0}$ and consider $\bar
x':=(2\lambda_0-\bar y, \bar
z)$,
the symmetric of $\bar x$ with respect to the hyperplane $T_{\lambda_0}$. By
reflecting $\bar x'$ through the hyperplanes $T^{\nu}_{\lambda_0(\nu)}$
for $\nu\in I_{\delta_2}(e_1)$, we obtain the points
$$
A(\nu)=(y(\nu),z(\nu))=\bar x'+2(\lambda_0(\nu)-\bar x'.\nu)\nu,
$$
which, for what we remarked before, belong to $\partial C$ . Indeed,
$\bar x\in \partial C\cap\Omega_{\lambda_0}$, so $\bar x\in \partial
C'$, $\bar x '\in\partial
(R_{\lambda_0}(C'))=\partial(R_{\lambda_0(\nu)}^\nu(C'))
$ and so $A(\nu)\in\partial C'$.
Since $\bar x\in \Omega_{\lambda_0}$, we can suppose, taking $\delta_2$
smaller if necessary, that for each $\nu\in I_{\delta_2}(e_1)$,
the point $A(\nu)$ belongs to $\Omega_{\lambda_0}$ and that
$\lambda_0(\nu)-\bar
x'.\nu<0$ (since $\bar x'.e_1>\lambda_0$ and $\nu\mapsto\lambda_0(\nu)$ is
continuous). So $A(\nu)\in\partial C$ for $\nu\in I_{\delta_2}(e_1)$ (see
the definition of
$C$ and $C'$). Observe that
by the continuity of $\nu\mapsto\lambda_0(\nu)$,
the function $\nu\mapsto A(\nu)$ is also continuous
and it is injective as it is easy to see.

By (\ref{4.8}), $u(\bar x)=u(\bar x')=u(A(\nu))$ for each $\nu\in
I_{\delta_2}(e_1)$, so that the function $u$ is constant on the set
$\Gamma:=\{A(\nu)\,|
\,\nu\in I_{\delta_2}(e_1)\}$. Since $A(\nu)\in \partial C\cap
\Omega_{\lambda_0}\subset\partial C'$
for all $\nu\in I_{\delta_2}(e_1)$, we have by (\ref{propr}) $\nabla u=0$ on
$\Gamma$. We will prove that the projection of $\Gamma$ on the
hyperplane $T_{\lambda_0}$ contains an open subset of
$T_{\lambda_0}$, this will lead to a contradiction with Lemma \ref{prop 3.1}
and the proof of Theorem \ref{mono} will be concluded.

Let us now write the generic direction $\nu\in S^{N-1}:=\{v\in\R^N\,|\,|v|=1\}$
as $\nu=(\nu_y,\nu_z)$, with $\nu_y\in\R$, $\nu_z\in\R^{N-1}$. If $\nu$
is close to $e_1$, then $\nu_y=\sqrt{1-|\nu_z|^2}$.

We take now $\beta>0$ small and  consider the set
$$K:=\left\{\nu=(\nu_y,
\nu_z)\,|\,\nu_z\in\bar B_\beta,\nu_y=\sqrt{1-|\nu_z|^2}\right\},$$
where $\bar B_\beta=\left\{z\in\R^{N-1}\,|\,|z|\le\beta\right\}$
is the closed ball in $\R^{N-1}$ centered at the origin with
radius $\beta$.

By construction, $K$ is a compact neighborhood of $e_1$ in the
metric space $S^{N-1}$, and if $\beta$ is small enough, then $K$
is contained in $I_{\delta_2}(e_1)$.
We will show that if $A(\nu)=(y(\nu),z(\nu))$ then the set
$\{z(\nu)\,|\,\nu\in K\}$
contains an open set in $\R^{N-1}$.

By definition of $A(\nu)$, $z(\nu)=\bar z'+2(\lambda_0(\nu)-\bar
x'.\nu)\nu_z$, where $\nu=(\sqrt{1-|\nu_z|^2},\nu_z)\in K$,
$\nu_z\in\bar B_\beta$. We will prove that the image of the function
$$
F(\nu_z):=2(\lambda_0(\nu)-\bar
x'.\nu)\nu_z,\quad\nu_z\in\bar B_\beta,\quad\nu=(\sqrt{1-|\nu_z|^2},\nu_z)
$$
contains a $(N-1)$-dimensional ball centered at the origin. This
will imply that $\{z(\nu)\,|\,\nu\in K\}$ contains an open set.

Let us consider a point $l\in S^{N-1}:=\{z\in\R^{N-1}\,|\,|z|=1\}$
and the segment $S_l:=\{tl\,|\,|t|\le\beta\}\subset\bar B_\beta$.
By definition of $F$ and since $F$ is continuous, the image $F(S_l)$
is a segment contained in the line passing
through the origin with direction $l$ in $R^{N-1}$. Moreover,
since $\lambda_0(\nu)-\bar x'.\nu<0$ for all $\nu\in K$ and since
$S_l$ contains points $\nu_z=tl$ with $t$
both positive and negative, the origin is an interior point of
$F(S_l)$. Hence we can write
$$
F(S_l)=\{tl\,|\,t\in[d_1(l),d_2(l)]\}\qquad d_1(l)<0<d_2(l),$$
for all $l\in S^{N-1}$ and for some functions $d_1(l),d_2(l)$.
 By the compactness of $K$ and the continuity of $\lambda_0(\nu)$ with
 respect to $\nu$, there exists $\delta\in(0,1)$ such that
$\lambda_0(\nu)-\bar x'.\nu<-\delta<0$ for all $\nu\in K$,
so that $d_1(l)\le -\beta\delta$ and $d_2(l)\ge \beta\delta$
for all $l\in S^{N-2}$. Thus 
the set
$$\{z\in\R^{N-1}\,|\,|z|\le \beta\delta\}=\bar B_{\beta\delta}\subseteq
F(\bar B_\beta),$$
which ends the proof.
 \cqfd
\mbox{}\\

A direct consequence of Theorem \ref{mono} is the following
symmetry result (Theorem \ref{corol} of the introduction).
\begin{cor}
Let $\nu\in\R^N$ and $\Omega\subset\R^N$ be a
domain with $C^1$ boundary symmetric
with respect to the hyperplane $T_0^\nu=\{x\in\R^N\,|\,x.\nu=0\}$
and $\lambda_1(\nu)=\lambda_1(-\nu)=0$.
Assume that one of the following conditions holds
\begin{enumerate}
\item $p_1,p_2\in(1,2)$ and  $f,g:\R\to\R^+$ are  strictly
increasing functions on $\R^+$ such that $f(x)>0, g(x)>0$ for all $x>0$,
\item $p_1\in(1,\infty),p_2=2$ and $f,g:\R\to\R^+$ are nondecreasing
on $\R^+$.
\end{enumerate}
Moreover suppose that $f$ and $g$ are locally Lipschitz continuous
on $\R$. Let $(u,v)  \in C^1_0(\Om)\times C^1_0(\Om)$ is
a weak solution of
(\ref{1.27}), then $u$ and $v$ are symmetric and decreasing in the $\nu$
direction in $\Omega_0^\nu$.  In
particular, if $\Omega$ is the ball $B_R(0)$ in $\R^N$ with center at the
origin
and radius $R$,
then $u,v$ are radially symmetric. If moreover $f(x)>0, g(x)>0$ for all $x>0$,
then $u'(r),v'(r)<0$ for $0<r<R$, $r=|x|$.
\end{cor}
\proof The symmetry of $u,v$ is a direct consequence of Theorems
\ref{5.3} and \ref{mono}. Suppose that $\Omega$ is a ball. If $r\in(0,R)$
and $G:=B_R\setminus \bar
B_r$, then $m:=u(r)$ is the maximum of $u$ in $\bar G$ and the
minimum of $u$ in $\bar B_r$. Since $f(x)>0,g(x)>0$
for all $x>0$,
$$
-\plun (u-m)=-\plun u>0\quad\textrm{in }B_r,
$$
and by the Hopf's lemma, $u'(r)<0$. A similar result shows that $v'(r)<0$
in $B_r$.\cqfd

\end{document}